\documentclass[12pt]{article}
\usepackage{amsfonts,amsmath,amssymb,epsf}
\usepackage{hyperref}
\usepackage{graphics} 
\usepackage{epstopdf} 
\usepackage{xypic}

\advance\textheight by \topskip
\newcounter{thm}

\newcommand\sect[1]{\section{#1}\setcounter{equation}0\setcounter{thm}0} 
                   
\newtheorem{thm}{Theorem}[section]
\newtheorem{defn}[thm]{Definition}
\newtheorem{prop}[thm]{Proposition}

\newtheorem{rema}[thm]{Remark}
\newtheorem{lemma}[thm]{Lemma}

\newcommand\void[1]       {}

\newcommand\be            {\begin{equation}}
\newcommand\bea           {\begin{equation}\begin{array}l\displaystyle}
\newcommand\bearll        {\begin{array}{ll}\displaystyle}
\newcommand\ee            {\end{equation}}
\newcommand\eear          {\end{array}}
\newcommand\enl           {\\[1em]\displaystyle}
\newcommand\etb           {&\!\! \displaystyle}
\newcommand\erf[1] {(\ref{#1})}
\newcommand\labl[1]       {\label{#1}\ee}

\renewcommand\cir         {\,{\circ}\,}
\newcommand\CxC           {\Cc\,{\boxtimes}\,\tilde\Cc}
\newcommand\eps           {\varepsilon}
\newcommand\Hom           {\mathrm{Hom}}
\newcommand\id            {{\rm id}}
\newcommand\In            {\,{\in}\,}

\newcommand\one           {{\bf1}}
\newcommand\oti           {\,{\otimes}\,}
\newcommand\ti            {\,{\times}\,}

\newcommand{\halmos}{\rule{1ex}{1.4ex}}
\newcommand{\pf}{\noindent{\it Proof.}\hspace{2ex}}
\newcommand{\epf}{\hspace*{\fill}\mbox{$\halmos$}}

\newcommand\Cb            {\mathbb{C}}
\newcommand\Nb            {\mathbb{N}}

\newcommand\Cc            {\mathcal{C}}
\newcommand\Ic            {\mathcal{I}}
\newcommand\Jc            {\mathcal{J}}

\begin{document}
\thispagestyle{empty}
\def\thefootnote{\fnsymbol{footnote}}
\begin{flushright}
KCL-MTH-07-10\\
0708.1897 [math.CT]
\end{flushright}
\vskip 1.0em
\begin{center}\LARGE
Morita classes of algebras in modular tensor categories
\end{center}\vskip 1.5em
\begin{center}\large
  Liang Kong%
  $^{a}$\footnote{Email: {\tt kong@mpim-bonn.mpg.de}}
  and
  Ingo Runkel%
  $^{b}$\footnote{Email: {\tt ingo.runkel@kcl.ac.uk}}%
\end{center}
\begin{center}\it$^a$
Max-Planck-Institut f\"ur Mathematik \\
Vivatsgasse 7, 53111 Bonn, Germany
\end{center}
\begin{center}\it$^b$
  Department of Mathematics, King's College London \\
  Strand, London WC2R 2LS, United Kingdom  
\end{center}
\vskip .5em
\begin{center}
  August 2007
\end{center}
\vskip 1em
\begin{abstract}
We consider algebras in a modular tensor category $\Cc$. If the 
trace pairing of an algebra $A$ in $\Cc$ is non-degenerate we
associate to $A$ a commutative algebra $Z(A)$, called the full
centre, in a doubled version of the category $\Cc$. We prove that two simple algebras with non-degenerate trace pairing are Morita-equivalent if and only if their full centres are isomorphic as algebras. This result has an interesting interpretation in two-dimensional rational conformal field theory; it implies that there cannot be several incompatible sets of boundary conditions for a given bulk theory.
\end{abstract}

\setcounter{footnote}{0}
\def\thefootnote{\arabic{footnote}}

\newpage

\sect{Introduction and summary}\label{sec:intro}

It is well-known that two Morita equivalent rings have isomorphic centres (see e.g.\ \cite[\S\,21]{and-ful}). The converse is in general not true, a counter example is provided by the real numbers and the quaternions. On the other hand, for simple algebras over $\Cb$ (or any algebraically closed field) the converse holds trivially, since all such algebras are of the form $\text{Mat}_n(\Cb)$ and all have centre $\Cb$. 

The situation becomes much richer if instead of considering algebras only in the category of finite-dimensional $\Cb$-vector spaces one allows for more general tensor categories. For example, for the categories of integrable highest weight representations of the affine Lie algebras $\widehat{\mathrm{sl}}(2)_k$, $k=1,2,\dots$, one finds an ADE-pattern for the Morita-classes, see e.g.\ \cite{ost}. These representation categories are in fact examples of so-called modular tensor categories, which are the class of categories we are considering in this paper.

We call an algebra non-degenerate if the trace pairing (or rather the appropriate categorical formulation thereof) is non-degenerate. We prove in this paper that two simple non-degenerate algebras in a modular tensor category are Morita equivalent if and only if they have isomorphic `full centres'. The latter is a commutative algebra which is a generalisation of the centre of an algebra over $\Cb$, but which typically lives in a different category than the algebra itself.

Our motivation to study the relation between Morita classes of algebras and their centres comes from two-dimensional conformal field theory. It has recently become clear that there is a close relationship between rational CFT and non-degenerate algebras in modular tensor categories, both in the Euclidean and the Minkowski formulation of CFT, see e.g.\ \cite{klm,tft1,lr,hu2,ko-cardy,unique}. In the Euclidean setting, the modular tensor category arises as the category of representations of a vertex operator algebra with certain additional properties \cite{hu1,hu2}, which we will refer to as `rational'. The non-degenerate algebra $A$ then is an algebra of boundary fields \cite{tft1}, i.e.\ an open-string vertex operator algebra \cite{hk1}. It turns out that $A$ and the rational vertex operator algebra together uniquely determine a CFT \cite{tft1,tft5,unique}; however, to ensure its existence, some complex analytic and convergence issues remain to be settled. As a consequence of the uniqueness, one can obtain from $A$ the algebra of bulk fields \cite{unique}, i.e.\ a full field algebra \cite{hk2}. An important question then is if two non-Morita equivalent open-string vertex operator algebras can give rise to the same full field algebra, or -- in more physical terms -- if there may exist several incompatible sets of boundary conditions for a given bulk CFT. Our result implies that for a CFT which is rational (in the sense that its underlying vertex operator algebra is rational), this cannot happen.

\bigskip

Recall that an algebra in a tensor category $\Cc$ with associator $\alpha_{U,V,W}$ and unit constraints $l_U,r_U$ is a triple $A=(A,m,\eta)$ where $A$ is an object of $\Cc$, $m$ (the multiplication) is a morphism $A \oti A \rightarrow A$ such that $m \cir (m \oti \id_A) \cir \alpha_{A,A,A} = m \cir (\id_A \oti m)$, and $\eta$ (the unit) is a morphism $\one \rightarrow A$ such that $m \cir (\id_A \oti \eta) = \id_A  \cir r_A$ and $m \cir (\eta\oti \id_A) = \id_A \cir l_A$. We will only consider unital algebras. In the following we will also assume that all tensor categories are strict to avoid spelling out associators and unit constraints.

In the same way one defines left--, right--, and bimodules. For example, given two algebras $A$ and $B$, an $A$-$B$-bimodule is a triple $X=(X,\rho_l,\rho_r)$ where $\rho_l : A \oti X \rightarrow X$ and $\rho_r : X \oti B \rightarrow X$ are the representation morphisms; 
$\rho_l$ obeys $\rho_l \cir (m_A \oti \id_X) = \rho_l \cir (\id_A \oti \rho_l)$ and $\rho_l \cir (\eta_A \oti \id_X) = \id_X$, and similar for $\rho_r$. Furthermore the left and right action commute, i.e.\ $\rho_r \cir (\rho_l \oti \id_B) = \rho_l \cir (\id_A \oti \rho_r)$.

With the help of bimodules we can now define when an algebra is simple, namely when it is simple as a bimodule over itself, and when two algebras $A$, $B$ are Morita equivalent, namely when there exist an $A$-$B$-bimodule $X$ and a $B$-$A$-bimodule $Y$ such that $X \otimes_B Y \cong A$ and $Y \otimes_A X \cong B$ as bimodules.

\medskip

Let now $\Cc$ be a modular tensor category (see \cite{tur} and e.g.\ \cite{baki}), i.e.\ a semisimple $\Cb$-linear abelian ribbon category with $\mathrm{End}(\one) = \Cb\,\id_\one$, having a finite number of isomorphism classes of simple objects and whose braiding obeys a certain nondegeneracy condition. (This definition is slightly more restrictive than the original one in \cite{tur}.) We will express morphisms in ribbon categories with the help of the usual graphical notation \cite{js}; our conventions are summarised in \cite[app.\,A.1]{tft5}. Given an algebra $A$ we can define the morphism $\Phi_A : A \rightarrow A^\vee$ as
\be
  \Phi_A ~=~ 
  \raisebox{-52pt}{
  \begin{picture}(80,108)
   \put(0,8){\scalebox{.75}{\includegraphics{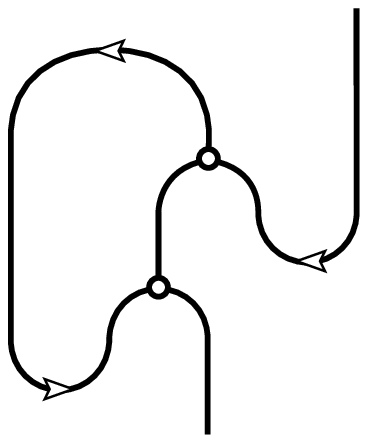}}}
   \put(0,8){
     \setlength{\unitlength}{.75pt}\put(-151,-235){
     \put(215,320)  {\scriptsize $ m $}
     \put(200,283)  {\scriptsize $ m $}
     \put(205,224)  {\scriptsize $ A $}
     \put(183,293)  {\scriptsize $ A $}
     \put(207,340)  {\scriptsize $ A $}
     \put(248,363)  {\scriptsize $ A^\vee $}
     }\setlength{\unitlength}{1pt}}
  \end{picture}}
  ~~~~.
\labl{eq:PhiA-def}
As in \cite{bim} we call an algebra $A$ {\em non-degenerate} iff $\Phi_A$ is invertible (the definition still makes sense in a tensor category with dualities). This generalises the condition that the trace pairing $a,b \mapsto \mathrm{tr}(a \cdot b)$ of a finite-dimensional algebra over a field is non-degenerate. We will list some properties of non-degenerate algebras in section \ref{sec:pre-nondeg} below.

Given an algebra $A$ in $\Cc$, the non-trivial braiding leads to two notions of centre, namely the left centre $C_l(A)$ and the right centre $C_r(A)$ of $A$
\cite{vz,ost,corr}. Denoting the braiding of $\Cc$ by $c_{U,V} : U \oti V \rightarrow V \oti U$, the left centre is the largest subobject
$C_l(A) \overset{\iota_l}{\rightarrow} A$ such that $m \circ c_{A,A} \circ (\iota_l \oti \id_A) = m \circ (\iota_l \oti \id_A)$ and the right centre the largest subobject $C_r(A) \overset{\iota_r}{\rightarrow} A$ such that $m \circ c_{A,A} \circ (\id_A \oti \iota_r) = m \circ (\id_A \oti \iota_r)$. We will give a formulation of the left centre of a non-degenerate algebra as the image of an idempotent in section \ref{sec:pre-full} below.

\medskip

The final ingredient we need to state our main result is a doubled version of $\Cc$, namely $\CxC$. Here the category $\tilde\Cc$ is obtained from $\Cc$ by replacing braiding and twist with their inverses, and the product $\CxC$ is the completion with respect to direct sums of $\Cc\ti\tilde\Cc$ 
(where the objects are pairs of objects in $\Cc$ and the $\Hom$-spaces are tensor products of the two corresponding $\Hom$-spaces in $\Cc$). $\CxC$ is again a modular tensor category. 
In fact, there is a notion of a `centre' $\mathcal{Z}$ of a tensor category, and for a modular tensor category $\Cc$ one finds $\mathcal{Z}(\Cc) \cong \CxC$ \cite{mug2}.

Apart from the tensor unit, the category $\CxC$ contains another canonically given commutative non-degenerate algebra, defined on the object $R = \oplus_{i\in\Ic} \, U_i \ti U^\vee_i$ \cite{mug2,corr,ko-ffa}. Here the (finite) set $\Ic$ indexes a choice of representatives $U_i$ of the isomorphism classes of simple objects in $\Cc$. The multiplication and further properties of $R$ are given in section \ref{sec:pre-mod}.

For a non-degenerate algebra $A$ in $\Cc$ we can now define the {\em full centre} $Z(A)$ as the left centre of the algebra $(A \ti \one) \oti R$ in $\CxC$ \cite{tftcft,unique}; our convention for the tensor product of algebras and some properties of the full centre will be discussed in section \ref{sec:pre-full}. As opposed to the left and right centres, the full centre is not a subobject of $A$, in fact it is not even an object of the same category. On the other hand, one can recover $C_l(A)$ and $C_r(A)$ from $Z(A)$ by applying suitable projections. Furthermore, if $\Cc$ is the category $\mathcal{V}ect_f(\Cb)$ of finite-dimensional complex vector spaces then also $\CxC \cong \mathcal{V}ect_f(\Cb)$, and the notions of left, right and full centre coincide and agree with the usual definition of the centre of an algebra over a field.

The full centre turns out to be a Morita-invariant notion and our main result is that it can be used to distinguish Morita-classes of non-degenerate algebras. 

\begin{thm}\label{thm:main}
Let $\Cc$ be a modular tensor category and let $A$, $B$ be simple non-degenerate algebras in $\Cc$. Then the following two statements are equivalent.
\\
(i)\phantom{i} $A$ and $B$ are Morita equivalent.
\\
(ii) $Z(A)$ and $Z(B)$ are isomorphic as algebras.
\end{thm}

\begin{rema}\label{rem:intro}\end{rema}
\vspace*{-.9em}
(i) In the special case $\Cc = \mathcal{V}ect_f(\Cb)$ a simple non-degenerate algebra is isomorphic to the full matrix algebra $\mathrm{Mat}_n(\Cb)$ for some $n$, and the full centre $Z$ is just the usual centre, which in the case of $\mathrm{Mat}_n(\Cb)$ is $\Cb$. The above theorem then just states that any two full matrix algebras over $\Cb$ are Morita equivalent. 
\\[.3em]
(ii) An algebra is called haploid iff $\dim\Hom(\one,A)=1$ \cite{fs-cat}. Denote by $C_\mathrm{max}(\CxC)$ the set of isomorphism classes $[B]$ of haploid commutative non-degenerate algebras $B$ in $\CxC$ which obey in addition $\dim(B) = \mathrm{Dim}(\Cc)$, where $\mathrm{Dim}(\Cc) = \sum_{i\in\Ic} \dim(U_i)^2$. (It follows from \cite[thm.\,4.5]{kios} that this is the maximal dimension such an algebra can have.) Note that $[R] \In C_\mathrm{max}(\CxC)$, with $R$ defined as above. Let further $M_\mathrm{simp}(\Cc)$ be the set of Morita classes $\{A\}$ of simple non-degenerate algebras $A$ in $\Cc$. We will see in remark \ref{rem:Z-prop}\,(ii) that the assignment $z : \{A\} \mapsto [Z(A)]$ is a well-defined map from $M_\mathrm{simp}(\Cc)$ to $C_\mathrm{max}(\CxC)$. For example, $z(\{\one\}) = [R]$. Theorem \ref{thm:main} shows that $z$ is injective. A result recently announced by M\"uger \cite{mug-conf} shows that $z$ is also surjective.
(An independent proof of surjectivity has subsequently appeared
in \cite[sect.\,3.3]{us2}.)
\\[.3em]
(iii) A closed two-dimensional topological field theory is the same as a commutative Frobenius algebra $B$ over $\Cb$, see e.g.\ \cite{kock}. In the case that $B$ is semi-simple, the possible boundary conditions for the theory defined by $B$ can be classified by $K_0( B\text{-mod} )$ \cite{mo,mose}. For a (rational) two-dimensional {\em conformal} field theory the boundary conditions can be classified by $K_0(A\text{-mod})$ where $A$ is a non-degenerate algebra in $\Cc$, and $\Cc$ in turn is the representation category of a rational vertex algebra $\mathcal{V}$ \cite{tft1}. The algebra $A$ comes from the boundary fields -- i.e.\ from an open-string vertex algebra over $\mathcal{V}$ -- for one of the possible boundary conditions \cite{tft1,hk1,ko-ocfa}. For the topological theory, the category $\Cc$ is given by $\Cc = \mathcal{V}ect_f(\Cb)$ and for $B$ one can choose the centre of $A$. (If $A$ is not simple this choice is not unique, see \cite{lp} and \cite[rem.\,4.27]{unique}.) For $\mathcal{V}ect_f(\Cb)$, $A$ and $B\,{=}\,Z(A)$ are Morita-equivalent, and so $K_0$ of $A$-mod and $B$-mod agree. In general one finds that, for $A$ a simple non-degenerate algebra in a modular tensor category $\Cc$ and $B=Z(A)$ the full centre,
$$
 \#\big(\text{isocl.\ of simple $B$-left modules in $\CxC$}\big)
 ~=~
 \#\big(\text{isocl.\ of simple $A$-$A$-bimodules $\Cc$}\big) ~.
$$
This can be computed from \cite[thm.\,5.18]{tft1} together with the fact that $Z(A)$ has a unique (up to isomorphism) simple local left module, namely $Z(A)$ itself. Thus in general, $K_0(B\text{-mod})$ -- the Grothendieck group of the category of $B$-left modules in $\CxC$ -- is related to defect lines (see \cite[rem.\,5.19]{tft1} and \cite{defect}), and its relevance for the classification of boundary conditions is special to the topological case. Nonetheless, there is a connection between $B$ and boundary conditions: We will see in section \ref{sec:ii-to-i} that via the tensor functor $T: \CxC \rightarrow \Cc$ one obtains an algebra $T(B)$ in $\Cc$ which is a direct sum of simple non-degenerate algebras, all of which are Morita-equivalent to $A$. In fact (cf.\ prop.\,\ref{prop:CA=TA} below) one has that $K_0(T(B)\text{-mod}) \cong K_0\big(A\text{-mod}\big)^{\!\times n}$, where $n$ is the number of isomorphism classes of simple $A$-left modules in $\Cc$.

\medskip

The rest of the paper is organised as follows. In section \ref{sec:prelim} we collect some results on non-degenerate algebras and the full centre. Section \ref{sec:i-to-ii} we prove that statement (i) in theorem \ref{thm:main} implies (ii) and in section \ref{sec:ii-to-i} we prove the converse.

\sect{Preliminaries}\label{sec:prelim}

\subsection{Properties of non-degenerate algebras}\label{sec:pre-nondeg}

Not all the properties discussed in this section require us to work with the full structure of a modular tensor category and we therefore state them in the appropriate context. However, all these properties do in particular hold for modular tensor categories.

\medskip

Let $\Cc$ be a (strict) tensor category. In the same way that one defines an algebra in $\Cc$ one can define a coalgebra $A = (A,\Delta,\eps)$ where $\Delta : A \rightarrow A \oti A$ and $\eps : A \rightarrow \one$ obey co-associativity and the counit condition.
\begin{defn} {\rm 
A Frobenius algebra $A = (A,m,\eta,\Delta,\eps)$ is an algebra and a coalgebra such that the coproduct is an intertwiner of $A$-bimodules, i.e.\ $(\id_A \oti m) \cir (\Delta \oti \id_A) = \Delta \oti m = (m \oti \id_A) \cir (\id_A \oti \Delta)$. 
}
\end{defn}

We will use the following graphical representation for the morphisms of a Frobenius algebra,
\be
  m = \raisebox{-20pt}{
  \begin{picture}(30,45)
   \put(0,6){\scalebox{.75}{\includegraphics{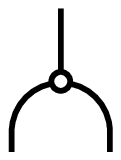}}}
   \put(0,6){
     \setlength{\unitlength}{.75pt}\put(-146,-155){
     \put(143,145)  {\scriptsize $ A $}
     \put(169,145)  {\scriptsize $ A $}
     \put(157,202)  {\scriptsize $ A $}
     }\setlength{\unitlength}{1pt}}
  \end{picture}}  
  ~~,\quad
  \eta = \raisebox{-15pt}{
  \begin{picture}(10,30)
   \put(0,6){\scalebox{.75}{\includegraphics{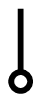}}}
   \put(0,6){
     \setlength{\unitlength}{.75pt}\put(-146,-155){
     \put(146,185)  {\scriptsize $ A $}
     }\setlength{\unitlength}{1pt}}
  \end{picture}}
  ~~,\quad
  \Delta = \raisebox{-20pt}{
  \begin{picture}(30,45)
   \put(0,6){\scalebox{.75}{\includegraphics{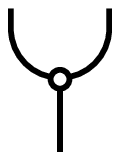}}}
   \put(0,6){
     \setlength{\unitlength}{.75pt}\put(-146,-155){
     \put(143,202)  {\scriptsize $ A $}
     \put(169,202)  {\scriptsize $ A $}
     \put(157,145)  {\scriptsize $ A $}
     }\setlength{\unitlength}{1pt}}
  \end{picture}}
  ~~,\quad
  \eps = \raisebox{-15pt}{
  \begin{picture}(10,30)
   \put(0,10){\scalebox{.75}{\includegraphics{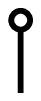}}}
   \put(0,10){
     \setlength{\unitlength}{.75pt}\put(-146,-155){
     \put(146,145)  {\scriptsize $ A $}
     }\setlength{\unitlength}{1pt}}
  \end{picture}}
  ~~.
\ee
A Frobenius algebra $A$ in a $\Bbbk$-linear tensor category, for a field $\Bbbk$, is called {\em special} iff $m \circ \Delta = \zeta \, \id_A$ and $\eps \circ \eta = \xi \, \id_\one$ for nonzero constants $\zeta$, $\xi \In \Bbbk$. If $\zeta=1$ we call $A$ {\em normalised-special}. 

A (strictly) sovereign tensor category is a tensor category equipped with a left and a right duality which agrees on objects and morphisms (see e.g.\ \cite{bichon,fs-cat} for more details). We will write the dualities as
\be
\begin{array}{llll}
  \raisebox{-8pt}{
  \begin{picture}(26,22)
   \put(0,6){\scalebox{.75}{\includegraphics{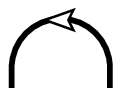}}}
   \put(0,6){
     \setlength{\unitlength}{.75pt}\put(-146,-155){
     \put(143,145)  {\scriptsize $ U^\vee $}
     \put(169,145)  {\scriptsize $ U $}
     }\setlength{\unitlength}{1pt}}
  \end{picture}}  
  \etb= d_U : U^\vee \oti U \rightarrow \one
  ~~,\qquad &
  \raisebox{-8pt}{
  \begin{picture}(26,22)
   \put(0,6){\scalebox{.75}{\includegraphics{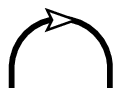}}}
   \put(0,6){
     \setlength{\unitlength}{.75pt}\put(-146,-155){
     \put(143,145)  {\scriptsize $ U $}
     \put(169,145)  {\scriptsize $ U^\vee $}
     }\setlength{\unitlength}{1pt}}
  \end{picture}}  
  \etb= \tilde d_U : U \oti U^\vee \rightarrow \one
  ~~,
\\[2em]
  \raisebox{-8pt}{
  \begin{picture}(26,22)
   \put(0,0){\scalebox{.75}{\includegraphics{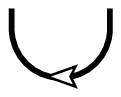}}}
   \put(0,0){
     \setlength{\unitlength}{.75pt}\put(-146,-155){
     \put(143,183)  {\scriptsize $ U $}
     \put(169,183)  {\scriptsize $ U^\vee $}
     }\setlength{\unitlength}{1pt}}
  \end{picture}}  
  \etb= b_U : \one \rightarrow U \oti U^\vee
  ~~,
  &
  \raisebox{-8pt}{
  \begin{picture}(26,22)
   \put(0,0){\scalebox{.75}{\includegraphics{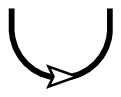}}}
   \put(0,0){
     \setlength{\unitlength}{.75pt}\put(-146,-155){
     \put(143,183)  {\scriptsize $ U^\vee $}
     \put(169,183)  {\scriptsize $ U $}
     }\setlength{\unitlength}{1pt}}
  \end{picture}}  
  \etb= \tilde b_U : \one \rightarrow U^\vee \oti U
  ~,
\eear\ee
(n.b.,\ `b' stands for birth and `d' for death).
Given these dualities one can define the left and right traces of a morphism $f: U \rightarrow U$ as 
$\mathrm{tr}_l(f) = d_U \cir (\id_{U^\vee} \oti f) \cir \tilde b_U$
and
$\mathrm{tr}_r(f) = \tilde d_U \cir (f \oti \id_{U^\vee}) \cir b_U$,
as well as the left and right dimension of $U$,
$\dim_{l/r}(U) = \mathrm{tr}_{l/r}(\id_U)$. 
If $U \cong U^\vee$, then $\dim_{l}(U)=\dim_{r}(U)$
\cite[rem.\,3.6.3]{fs-cat}.
In a modular tensor category
(and more generally in a spherical category) the left and right traces and dimensions always coincide.

Let now $\Cc$ be a sovereign tensor category.
A Frobenius algebra in $\Cc$ is {\em symmetric} iff
\be
  \raisebox{-35pt}{
  \begin{picture}(50,75)
   \put(0,8){\scalebox{.75}{\includegraphics{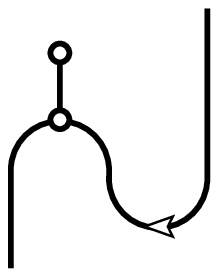}}}
   \put(0,8){
     \setlength{\unitlength}{.75pt}\put(-34,-37){
     \put(31, 28)  {\scriptsize $ A $}
     \put(87,117)  {\scriptsize $ A^\vee $}
     }\setlength{\unitlength}{1pt}}
  \end{picture}}
  ~=~
  \raisebox{-35pt}{
  \begin{picture}(50,75)
   \put(0,8){\scalebox{.75}{\includegraphics{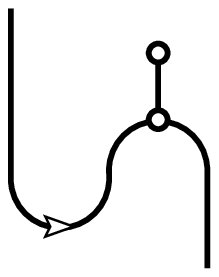}}}
   \put(0,8){
     \setlength{\unitlength}{.75pt}\put(-34,-37){
     \put(87, 28)  {\scriptsize $ A $}
     \put(31,117)  {\scriptsize $ A^\vee $}
     }\setlength{\unitlength}{1pt}}
  \end{picture}}
  ~~.
\labl{eq:Frob-sym-cond}
For a Frobenius algebra (symmetric or not) one has $\dim_l(A) = \dim_r(A)$ since the two morphisms in \erf{eq:Frob-sym-cond} are in fact isomorphisms (see e.g.\ the proof of lemma 3.7 in \cite{tft1}), and so $A \cong A^\vee$. We will write $\dim(A) \equiv \dim_{l/r}(A)$.

In a $\Bbbk$-linear sovereign category with $\mathrm{End}(\one) = \Bbbk \, \id_\one$ we will identify $\dim_{l/r}(U)$ with the corresponding element of $\Bbbk$ via $\mathrm{tr}_{l/r}(\id_U) = \dim_{l/r}(U) \, \id_\one$. In this case one finds that for a normalised-special symmetric Frobenius algebra one has $\eps \cir \eta = \dim(A) \, \id_\one$ \cite[sect.\,3]{fs-cat}; in particular, $\dim(A) \neq 0$.

\begin{defn}  {\rm 
An algebra $A$ in $\Cc$ is non-degenerate if 
the morphism $\Phi_A$ in \erf{eq:PhiA-def} is invertible.
}\end{defn}

The relation between non-degenerate algebras and Frobenius algebras is
summarised in the following lemma.

\begin{lemma}\label{lem:nondeg_vs_Frob}
Let $\Cc$ be a sovereign tensor category.
\\[.3em]
(i)\phantom{ii} Let $A$ be a non-degenerate algebra in $\Cc$. Taking
$\Delta = (\Phi_A^{-1} \oti m) \cir (\tilde b_A \oti \id_A)$ and $\eps = \eta^\vee \cir\Phi_A$ turns $A$ into a symmetric Frobenius algebra which obeys $m \cir \Delta = \id_A$ and $\eps \cir \eta = \dim(A)$.
\\[.3em]
(ii)\phantom{i} Let $A$ be a symmetric Frobenius algebra in $\Cc$ such that $m \cir \Delta = \id_A$. Then $A$ is a non-degenerate algebra.
\\[.3em]
(iii) Two non-degenerate algebras $A$ and $B$ are isomorphic as algebras if and only if they are isomorphic as Frobenius algebras (with counit and coproduct as given in (i)).
\\[.3em] 
(iv)\, If $\Cc$ is in addition 
$\Bbbk$-linear with $\mathrm{End}(\one) = \Bbbk \, \id_\one$
and $A$ is a non-degenerate algebra in $\Cc$ with 
$\dim(A) \neq 0$, then $A$ is special.
\end{lemma}

This lemma can be proved by combining and adapting lemmas 3.7, 3.11 and 3.12 of \cite{tft1}. Part of the proof of (i) involves showing that $\Phi_A$ in \erf{eq:PhiA-def} is also equal to the morphism obtained by `reflecting' the graph along a vertical axis (cf.\ \cite[eqn.\,(3.33)]{tft1}), and that equally $\Delta = (m \oti\Phi_A^{-1}) \cir (b_A \oti \id_A)$. In this sense, the Frobenius algebra structure on a non-degenerate algebra does not involve any arbitrary choices.

Whenever we will consider a non-degenerate algebra as a Frobenius algebra we mean the coproduct and counit given in part (i) of the above lemma.

\medskip

In the setting we will work with below, $\Cc$ is a modular tensor category and one can convince oneself that a simple non-degenerate algebra in $\Cc$ necessarily has $\dim(A) \neq 0$, cf.\ \cite[lem.\,2.6]{defect}. In particular, a simple non-degenerate algebra is then always also normalised-special symmetric Frobenius.

\medskip

Let now $\Cc$ be an abelian sovereign tensor category. Let $A$ be a non-degenerate algebra in $\Cc$ and let $M$ be a right $A$-module and $N$ be a left $A$-module. The tensor product $M \otimes_A N$ can be written as the image of the idempotent
\be
  P_{\otimes A} ~=~ 
  \raisebox{-34pt}{
  \begin{picture}(65,75)
   \put(8,8){\scalebox{.75}{\includegraphics{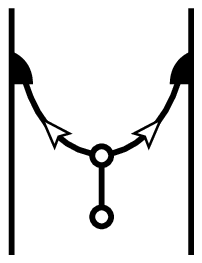}}}
   \put(8,8){
     \setlength{\unitlength}{.75pt}\put(-30,-38){
     \put(25, 28)  {\scriptsize $ M $}
     \put(25,115)  {\scriptsize $ M $}
     \put(79, 28)  {\scriptsize $ N $}
     \put(79,115)  {\scriptsize $ N $}
     \put(12, 87)  {\scriptsize $ \rho_M $}
     \put(86, 87)  {\scriptsize $ \rho_N $}
     \put(41, 62)  {\scriptsize $ A $}
     }\setlength{\unitlength}{1pt}}
  \end{picture}}
\ee
That is, there exist morphisms $e_A : M \otimes_A N \rightarrow M \oti N$ and $r_A : M \oti N \rightarrow M \otimes_A N$ such that $r_A \cir e_A = \id_{M \otimes_A N}$ and $e_A \cir r_A = P_{\otimes A}$. One can convince oneself that $r_A : M \oti N \rightarrow M \otimes_A N$ fulfils the universal property of the coequaliser of $\rho_M \oti \id_N$ and $\id_M \oti \rho_N$.

\subsection{Modular tensor categories}\label{sec:pre-mod}

Let $\Cc$ be a modular tensor category. Recall from section \ref{sec:intro} that we chose representatives $\{\, U_i \,|\, i \In \Ic \,\}$ for the isomorphism classes of simple objects. We also fix $U_0 = \one$ and for an index $k \in \Ic$ we define the index $\bar k$ by $U_{\bar k} \cong U^\vee_k$. The numbers $s_{i,j} \in \Cb$ are defined via
\be 
  s_{i,j} \, \id_\one = \mathrm{tr}(c_{U_i,U_j} \cir c_{U_j,U_i})\ . 
\ee
They obey $s_{i,j} = s_{j,i}$ and $s_{0,i} = \dim(U_i)$, see e.g.\ \cite[sect.\,3.1]{baki}. The non-degeneracy condition on the braiding of a modular tensor category is that the $|\Ic|{\times}|\Ic|$-matrix $s$ 
should be invertible. In fact, 
\be
  \sum_{k \in \Ic} s_{ik} \, s_{kj} = \mathrm{Dim}(\Cc) \, \delta_{i,\bar\jmath}
\ee 
(cf.\ \cite[thm.\,3.1.7]{baki}), where as above $\mathrm{Dim}(\Cc) = \sum_{i\in\Ic} \dim(U_i)^2$. In particular, $\mathrm{Dim}(\Cc) \neq 0$. One can show (even in the weaker context of fusion categories over $\Cb$) that $\mathrm{Dim}(\Cc) \ge 1$ \cite[thm.\,2.3]{eno}. 

Let us fix a basis 
$\{ \lambda_{(i,j)k}^{\alpha} \}_{\alpha=1}^{N_{ij}^k}$ 
in $\Hom (U_i\otimes U_j, U_k)$ and the dual basis
$\{ \Upsilon^{(i,j)k}_{\alpha} \}_{\alpha=1}^{N_{ij}^k}$ 
in $\Hom (U_k, U_i\otimes U_j)$. 
The duality of the bases means that 
$\lambda_{(i,j)k}^{\alpha} \cir \Upsilon^{(i,j)k}_{\beta}
 = \delta_{\alpha,\beta}\, \id_{U_k}$. We also fix
$\lambda_{(0,i)i} = \lambda_{(i,0)i} = \id_{U_i}$.
We denote the basis vectors graphically as follows: 
\be
\lambda_{(i,j)k}^{\alpha} = 
  \raisebox{-23pt}{
  \begin{picture}(30,52)
   \put(0,8){\scalebox{.75}{\includegraphics{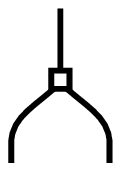}}}
   \put(0,8){
     \setlength{\unitlength}{.75pt}\put(-18,-11){
     \put(39,36)  {\scriptsize $ \alpha $}
     \put(30,61)  {\scriptsize $ U_k $}
     \put(15, 2)  {\scriptsize $ U_i $}
     \put(43, 2)  {\scriptsize $ U_j $}
     }\setlength{\unitlength}{1pt}}
  \end{picture}}
\quad , \qquad
\Upsilon_{\alpha}^{(i,j)k} =
  \raisebox{-23pt}{
  \begin{picture}(30,52)
   \put(0,8){\scalebox{.75}{\includegraphics{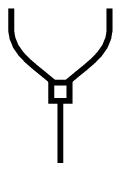}}}
   \put(0,8){
     \setlength{\unitlength}{.75pt}\put(-18,-11){
     \put(39,28)  {\scriptsize $ \alpha $}
     \put(30, 2)  {\scriptsize $ U_k $}
     \put(15,61)  {\scriptsize $ U_i $}
     \put(43,61)  {\scriptsize $ U_j $}
     }\setlength{\unitlength}{1pt}}
  \end{picture}} 
  ~~.
\ee
As in section \ref{sec:intro} let $R$ be the object in $\CxC$ given by 
$R = \oplus_{i \in \Ic}\, U_i \ti U_i^{\vee}$.
We define a unit morphism 
$\eta_R: \one \times \one \rightarrow R$ to be
the natural embedding and a multiplication morphism 
$m_R: R\otimes R \rightarrow R$ as
\be 
m_R ~=~ \bigoplus_{i,j,k\in \Ic} \sum_{\alpha=1}^{N_{ij}^{~k}} ~~
  \raisebox{-35pt}{
  \begin{picture}(120,80)
   \put(0,8){\scalebox{.75}{\includegraphics{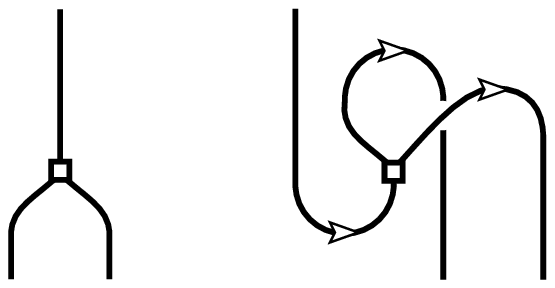}}}
   \put(0,8){
     \setlength{\unitlength}{.75pt}\put(-15,-15){
     \put(66,52)  {$ \times $}
     \put(36,48)  {\scriptsize $ \alpha $}
     \put(113,43)  {\scriptsize $ \alpha $}
     \put(27,98)  {\scriptsize $ U_k $}
     \put(11, 5)  {\scriptsize $ U_i $}
     \put(41, 5)  {\scriptsize $ U_j $}
     \put(93,98)  {\scriptsize $ U_k^\vee $}
     \put(135, 5)  {\scriptsize $ U_i^\vee $}
     \put(163, 5)  {\scriptsize $ U_j^\vee $}
     }\setlength{\unitlength}{1pt}}
  \end{picture}}
~~.  
\labl{eq:m-R}
The basis elements $\lambda_{(i,j)k}^{\alpha}$ and their duals enter the above
expression in such a way that $m_R$ is independent of the choice of
bases.
Using \erf{eq:m-R}, one can compute $\Phi_{R}$ 
(defined in \erf{eq:PhiA-def}) explicitly, resulting in
$\Phi_{R} = \mathrm{Dim}(\Cc) \bigoplus_{i\in \Ic} f_i \ti g_i$
where $f_i : U_i \rightarrow U_{\bar\imath}^\vee$ and
$g_i : U_i^\vee \rightarrow U_{\bar\imath}^{\vee\vee}$ are given by
\be
  f_i = (\lambda_{(i,\bar\imath)0}^{1} \oti \id_{U_{\bar\imath}^\vee})
  \cir (\id_{U_i} \oti b_{U_{\bar\imath}})
  ~~,\quad
  g_i = (\delta_{U_{\bar\imath}} \oti \tilde d_{U_i}) \cir
  \big( (c_{U_{\bar\imath},U_i}^{-1} \cir 
  \Upsilon_{1}^{(i,\bar\imath)0})
  \oti \id_{U_{i}^\vee} \big) ~,
\ee
and $\delta_U : U \rightarrow U^{\vee\vee}$ is the isomorphism
$(\tilde d_U \oti \id_{U^{\vee\vee}}) \cir (\id_U \oti b_{U^\vee})$.
It follows from \cite[prop.\,4.1]{mug2} 
(see also \cite[lem.\,6.19]{corr} and \cite[thm.\,5.2]{ko-ffa}) that 
the three morphisms $\eta_R$, $m_R$ and $\Phi_R$
give $R$ the structure of haploid commutative non-degenerate algebra. 
Thus it is also normalised-special symmetric Frobenius. (The algebra $R$ can also be defined in more general categories, see \cite{mug2,corr}.)

\subsection{Properties of the full centre}\label{sec:pre-full}

{}From hereon we will always take $\Cc$ to be a modular tensor category. Most of the constructions in this section can be carried out in greater generality, see e.g.\ \cite{corr}, but for the purpose of the proof of theorem \ref{thm:main} this will not be necessary.

\medskip

An algebra $A$ in a braided tensor category has  
a left centre and a right centre \cite{vz,ost}, both of which are sub-algebras of $A$. We will only need the left centre. The following definition is the one used in \cite{corr}, which in our setting is equivalent to that of \cite{vz,ost}. 

\begin{defn} {\rm 
The left centre $C_l(A)$ of a non-degenerate algebra $A$ in $\Cc$ is the image of the idempotent $P_l(A) : A \rightarrow A$, where
\be
  P_l(A) ~=~
  \raisebox{-40pt}{
  \begin{picture}(54,80)
   \put(0,8){\scalebox{.75}{\includegraphics{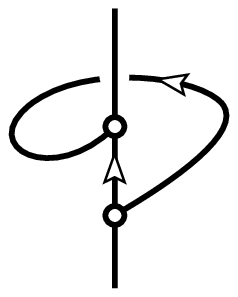}}}
   \put(0,8){
     \setlength{\unitlength}{.75pt}\put(-18,-19){
     \put(45, 10)  {\scriptsize $ A $}
     \put(45,105)  {\scriptsize $ A $}
     \put(78, 79)  {\scriptsize $ A $}
     \put(55, 65)  {\scriptsize $ m $}
     \put(36, 36)  {\scriptsize $ \Delta $}
     }\setlength{\unitlength}{1pt}}
  \end{picture}}
~~.
\labl{eq:Pl(A)-def}
}
\end{defn}

That is, there are morphisms $C_l(A) \overset{\iota_l}{\rightarrow} A$ 
and $A \overset{r_l}{\rightarrow} C_l(A)$ such that $r_l \cir \iota_l = \id_{C_l(A)}$ and $\iota_l \cir r_l = P_l(A)$. This follows from definition 2.31 and remark 2.34 of \cite{corr}. (In \cite{corr} the algebra $A$ is assumed to be special, but the relevant proofs only use $m \cir \Delta = \id_A$, which is satisfied by a non-degenerate algebra according to lemma \ref{lem:nondeg_vs_Frob}\,(i).) The proof of the following lemma can be found in \cite[prop.\,2.37]{corr}. 

\begin{lemma} \label{lem:leftcentre}
Let $A$ be a non-degenerate algebra in a modular tensor category $\Cc$.
\\[.3em]
(i)\phantom{i} $C_l(A)$ is a commutative symmetric Frobenius algebra.
\\[.3em]
(ii) If $C_l(A)$ is simple and $\dim(C_l(A)) \neq 0$, then $C_l(A)$ is in addition special.
\end{lemma}

\medskip

Given two algebras $A$ and $B$ we define a multiplication on the tensor product $A \oti B$ as $m_{A\otimes B} = (m_A \oti m_B) \cir (\id_A \oti c^{-1}_{A,B} \oti \id_B)$ and a unit morphism as $\eta_{A \otimes B} = \eta_A \otimes \eta_B$. This turns $A\oti B$ into an algebra. Note that one can also define a different multiplication $m_{A\otimes B}'$ by using $c_{B,A}$ instead of $c^{-1}_{A,B}$. The resulting algebra is isomorphic to $(A_\mathrm{op} \oti B_\mathrm{op})_\mathrm{op}$, where `op' stands for the opposed algebra, see \cite[rem.\,3.23]{tft1}. We will always use $m_{A\otimes B}$.

For two coalgebras we similarly set $\Delta_{A\otimes B} = (\id_A \oti c_{A,B} \oti \id_B) \cir (\Delta_A \oti \Delta_B)$ and $\eps_{A \otimes B} = \eps_A \oti \eps_B$. This turns $A \oti B$ into a coalgebra. One easily checks that if $A$ and $B$ share any of the properties non-degenerate, Frobenius, symmetric, special, then the property is inherited by $A \oti B$. On the other hand, even if $A$ and $B$ are commutative, $A \oti B$ is generally not.

\medskip

For an object $U$ of $\Cc$ denote by $R(U)$ the object in $\CxC$ given by $R(U) = (U \ti \one) \oti R$. 
($R(\,\cdot\,)$ can be understood as the adjoint of the functor
$T$ mentioned in remark \ref{rem:intro}\,(iii); more 
details can be found in \cite[sect.\,2.4]{us2}.)
If $A$ is a non-degenerate algebra in $\Cc$ then $A \ti \one$ is a non-degenerate algebra in $\CxC$ and the above discussion gives a non-degenerate algebra structure on $R(A)$. 

\begin{defn}  {\rm \cite[def.\,4.9]{unique}
The full centre $Z(A)$ of $A$ is defined to be $C_l(R(A))$.
}
\end{defn}

\begin{prop}\label{prop:Z(A)-prop}
Let $A$ be a non-degenerate algebra in a modular tensor category $\Cc$.
\\
(i)\phantom{i} $Z(A)$ is a commutative symmetric Frobenius algebra with $\dim(Z(A)) = d \cdot \mathrm{Dim}(\Cc)$ for some integer $d \ge 1$.
\\
(ii) If $A$ is simple then $Z(A)$ is a haploid commutative non-degenerate algebra with $\dim(Z(A))$ $=$ $\mathrm{Dim}(\Cc)$. Furthermore, $Z(A)$ is normalised-special.
\end{prop}

\noindent
\pf The first statement in
part (i) follows from lemma \ref{lem:leftcentre}\,(i) together with the above observation that $R(A)$ is a non-degenerate algebra in $\CxC$. For the
statement about the dimension, let $Z_{ij} = \dim \Hom(Z(A),U_i \ti U_j)$. By combining \cite[eqn.\,(A.3)]{tftcft} (note that in \cite{tftcft} $Z(A)$ has a different meaning, namely the object given in eqn.\,(3.9) there) with eqn.\,(5.65) and theorem~5.1 of \cite{tft1} it follows that $\sum_{k \in \Ic} Z_{ik} s_{kj} = \sum_{l \in \Ic} s_{il} Z_{lj}$. Using this we can compute
\be\bearll
  \dim(Z(A)) 
  \etb= \sum_{i,j} Z_{ij} \dim(U_i) \dim(U_j)
  = \sum_{i,j} s_{0i} \, Z_{ij} \, s_{j0}
  = \sum_{j,k} Z_{0k} \, s_{kj} \, s_{j0}
  \enl
  \etb= \sum_{k} Z_{0k}\, \delta_{k,0} \,\mathrm{Dim}(\Cc)
  = Z_{00} \, \mathrm{Dim}(\Cc)\ .
\eear\labl{eq:non-deg-center-dim}
It follows from the equalities (A.2) in \cite{tftcft} that $Z_{00} = \dim \Hom_{A|A}(A,A)$, where $\Hom_{A|A}(\,\cdot\,,\,\cdot\,)$ denotes the space of bimodule intertwiners. Since $\id_A$ is a bimodule intertwiner we have $Z_{00} \ge 1$.
\\[.3em] 
For (ii) note that in the present setting, $A$ is simple iff it is absolutely simple, i.e.\ iff $\Hom_{A|A}(A,A) = \Cb\,\id_A$, which is equivalent to $Z_{00}=1$. Therefore, $A$ is simple iff $Z(A)$ is haploid. Since by assumption in (ii), $A$ is simple, \erf{eq:non-deg-center-dim} holds with $Z_{00}=1$. Recall from above that $\mathrm{Dim}(\Cc) \ge 1$, so that altogether we see that $Z(A)$ is simple (since it is haploid) and has nonzero dimension. By lemma \ref{lem:leftcentre}\,(ii), $Z(A)$ is then also special. We can rescale the coproduct (and the counit) to make $Z(A)$ normalised-special and it then follows from lemma \ref{lem:nondeg_vs_Frob}\,(ii) that $Z(A)$ is non-degenerate.
\epf

\sect{Morita equivalence implies isomorphic full centre}\label{sec:i-to-ii}

Let $A$, $B$ be two non-degenerate algebras in a modular tensor category $\Cc$. Given an $A$-$B$-bimodule $X$ define the morphism $Q_X : R(B) \rightarrow R(A)$ as
\be
 Q_X = 
 \raisebox{-55pt}{
  \begin{picture}(65,112)
   \put(0,8){\scalebox{.75}{\includegraphics{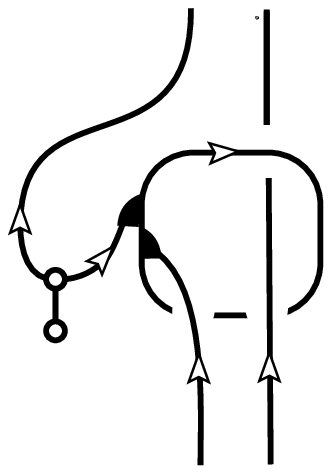}}}
   \put(0,8){
     \setlength{\unitlength}{.75pt}\put(-293,-160){
     \put(335,151)  {\scriptsize $ B {\times} \one$}
     \put(366,151)  {\scriptsize $ R $}
     \put(335,298)  {\scriptsize $ A {\times} \one$}
     \put(366,298)  {\scriptsize $ R $}
     \put(338,234)  {\scriptsize $ X {\times} \one$}
  }\setlength{\unitlength}{1pt}}
  \end{picture}}
\labl{eq:QX-def}
The morphism $Q_X$ is closely related to the linear map $D^{UV}_X$ defined in \cite{bim}, but is slightly more general as here we work with $A$-$B$-bimodules instead of $A$-$A$-bimodules.

\begin{lemma} \label{lem:QX-properties}
Let $A,B,C$ be non-degenerate algebras in $\Cc$, let $X,X'$ be $A$-$B$-bimodules and $Y$ a $B$-$C$-bimodule.
\\
(i)\phantom{ii} 
If $X \cong X'$ then $Q_X = Q_{X'}$.
\\
(ii)\phantom{i} 
$Q_A = P_l(R(A))$, with $P_l$ as defined in \erf{eq:Pl(A)-def}.
\\
(iii)
$Q_X \cir Q_Y = Q_{X \otimes_B Y}$.
\\
(iv)\phantom{i}  
$Q_X \cir P_l(R(B)) = Q_X = P_l(R(A)) \cir Q_X$.
\end{lemma}
\pf Part (i) is proved in the same way as the corresponding statement for $D^{UV}_X$, see \cite[eqn.\,(22)]{bim}. Namely, if $f : X \rightarrow X'$ is an isomorphism of bimodules, one inserts the identity $\id_X = f^{-1} \cir f$ anywhere on the $X$-loop in the pictorial representation \erf{eq:QX-def} of $Q_X$. One then drags $f$ around the loop until it combines with $f^{-1}$ to $f \cir f^{-1} = \id_{X'}$. This results in the morphism $Q_{X'}$.
\\
The equality in (ii) can be seen by comparing the pictorial representations and using that $R$ is (in particular) commutative and normalised-special; it also follows from the proof of \cite[prop.\,3.14(i)]{corr}. 
\\
Claim (iii) can be proved in the same way as \cite[lem.\,2]{bim}. Part (iv) is then a consequence of applying (i)--(iii) to $X \otimes_B B \cong X \cong A \otimes_A X$.
\epf

\medskip

Using $Q_X$ we define a morphism 
$D_X: Z(B) = C_l( (B\ti \one)\oti R) \rightarrow 
Z(A) = C_l( (A\ti \one)\oti R)$ 
by composing with the corresponding embedding and restriction morphisms,
\be
  D_X = r_l \cir Q_X \cir \iota_l \ .
\labl{eq:DX-def}
As a direct consequence of lemma \ref{lem:QX-properties} we have $D_X = D_{X'}$ for two isomorphic bimodules $X$ and $X'$, as well as, for $X$, $Y$ as in lemma \ref{lem:QX-properties},
\be
  D_A = \id_{Z(A)} 
  \quad , \qquad 
  D_X \circ D_Y = D_{X \otimes_B Y} \ .
\labl{eq:DX-properties}

\begin{lemma}\label{lem:phi_X-iso-frob}
Let $A, B$ be non-degenerate algebras (not necessarily simple)
and $X$ an $A$-$B$-bimodule, such that $\dim(A)$, $\dim(B)$ and $\dim(X)$ are
non-zero and the identities
\be 
 \raisebox{-28pt}{
  \begin{picture}(65,60)
   \put(0,8){\scalebox{.75}{\includegraphics{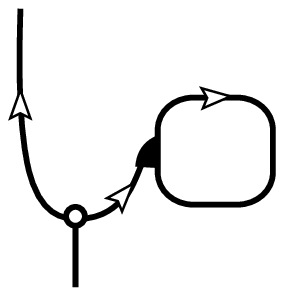}}}
   \put(0,8){
     \setlength{\unitlength}{.75pt}\put(-289,-193){
     \put(343,236)  {\scriptsize $ X$}
     \put(278,245)  {\scriptsize $ A$}
     \put(303,184)  {\scriptsize $ A$}
  }\setlength{\unitlength}{1pt}}
  \end{picture}}
 = \frac{\dim(X)}{\dim(A)} \,\, \id_A,  \hspace{1cm}
  \raisebox{-28pt}{
  \begin{picture}(65,60)
   \put(0,8){\scalebox{.75}{\includegraphics{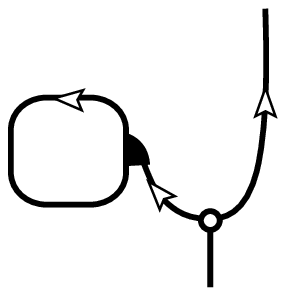}}}
   \put(0,8){
     \setlength{\unitlength}{.75pt}\put(-289,-193){
     \put(302,236)  {\scriptsize $ X$}
     \put(369,245)  {\scriptsize $ B$}
     \put(342,184)  {\scriptsize $ B $}
     }\setlength{\unitlength}{1pt}}
  \end{picture}}
 = \frac{\dim(X)}{\dim(B)} \,\, \id_B~,
\labl{eq:X-cond-unit}
and 
\be 
  \raisebox{-36pt}{
  \begin{picture}(60,80)
   \put(0,8){\scalebox{.75}{\includegraphics{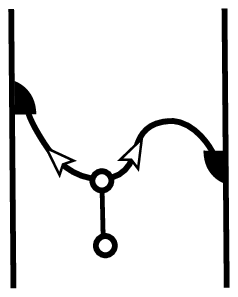}}}
   \put(0,8){
     \setlength{\unitlength}{.75pt}\put(-257,-293){
     \put(299,324)  {\scriptsize $ B $}
     \put(254,283)  {\scriptsize $ X $}
     \put(317,283)  {\scriptsize $ X^{\vee}$}
     \put(270,340)  {\scriptsize $ B $}
     }\setlength{\unitlength}{1pt}}
  \end{picture}} 
  =  \frac{\dim(A)}{\dim(X)} \,\,\, 
\raisebox{-42pt}{
\begin{picture}(100,90)
   \put(0,8){\scalebox{.75}{\includegraphics{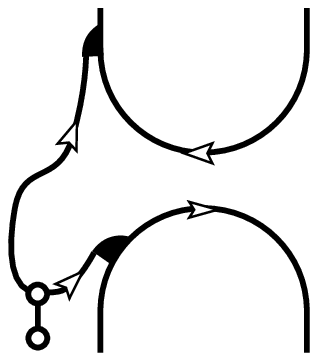}}}
   \put(0,8){
     \setlength{\unitlength}{.75pt}\put(-123,-162){
     \put(145,267)  {\scriptsize $ X $}
     \put(205,267)  {\scriptsize $ X^{\vee} $}
     \put(124,227)  {\scriptsize $ A $}
     \put(135,190)  {\scriptsize $ A $}
     \put(145,152)  {\scriptsize $ X $}
     \put(205,152)  {\scriptsize $ X^\vee $}
     }\setlength{\unitlength}{1pt}}
  \end{picture}}
\labl{eq:X-cond-mor}
hold. Then $\phi_X:= \frac{\dim(X)}{\dim(B)} D_X: Z(B)\rightarrow Z(A)$ 
is an isomorphism of Frobenius algebras. 
\end{lemma}

The precise form of the dimension-factors appearing in conditions \erf{eq:X-cond-unit} and \erf{eq:X-cond-mor} is not an extra condition, but is in fact uniquely fixed. For example composing the first equation in \erf{eq:X-cond-unit} with $\eps_A$ from the left and $\eta_A$ from the right gives the first constant.
Also note that $X^\vee$ is naturally a $B$-$A$-bimodule,
see e.g.\ \cite[sect.\,2.1]{defect}.

\medskip

\noindent{\it Proof of lemma \ref{lem:phi_X-iso-frob}}\\
a) {\it $A\cong X \otimes_B X^{\vee}$ as $A$-$A$-bimodules}: 
We define two morphisms $f_1: A \rightarrow X \otimes_B X^{\vee}$
and $f_2: X \otimes_B X^{\vee} \rightarrow A$ by 
\be\bearll
f_1 \etb= \frac{\dim(A)}{\dim(X)}\,\, 
r_B \cir (\rho_A \otimes \id_{X^{\vee}}) \cir  
(\id_A \otimes b_X)\ ,   \enl
f_2 \etb= (\id_A \otimes \tilde{d}_X) \cir 
(\id_A \otimes \rho_A \otimes \id_{X^{\vee}}) \cir 
((\Delta_A \cir \eta_A) \otimes e_B)\ .  
\eear\labl{f-1-def}
Notice first that both $f_1$ and $f_2$ are $A$-$A$-bimodule maps. 
It is easy to see that the condition
(\ref{eq:X-cond-unit}) implies $f_2\cir f_1 = \id_A$
and the condition (\ref{eq:X-cond-mor}) implies 
$f_1\cir f_2 = \id_{X \otimes_B X^{\vee}}$. Therefore,
$A\cong X \otimes_B X^{\vee}$ as bimodules and an isomorphism
is given by $f_1$. 
\\[.3em]
b) {\it $B\cong X^\vee \otimes_A X$ as $B$-$B$-bimodules}: 
This can be seen by a similar argument as used in a).
\\[.3em]
c) {\it $\phi_X$ is an isomorphism}: First note that taking the
trace of \erf{eq:X-cond-mor} and using \erf{eq:X-cond-unit} results
in the identity 
\be
  \dim(X)^2 = \dim(A) \dim(B)\ .
\labl{eq:dimX-dimAdimB}
Using this, as well as \erf{eq:DX-properties} and part b) we obtain
\be
  \phi_{X^\vee} \cir \phi_X
  = \frac{\dim(X)}{\dim(A)} \, \frac{\dim(X)}{\dim(B)} \, 
    D_{X^\vee} \cir D_X
  = D_{X^\vee \oti_A X} = D_B = \id_{Z(B)} ~,
\ee
In the same way one checks that $\phi_X \cir \phi_{X^\vee} = \id_{Z(A)}$.
Thus $\phi_X$ is an isomorphism.
\\[.3em]
d) {\it $\phi_X$ is an algebra map}:
The unit property $\phi_X \cir \eta_{Z(B)} = \eta_{Z(A)}$ can be seen as follows,
\be
 \frac{\dim(X)}{\dim(B)} 
 \raisebox{-62pt}{
  \begin{picture}(80,128)
   \put(0,8){\scalebox{.75}{\includegraphics{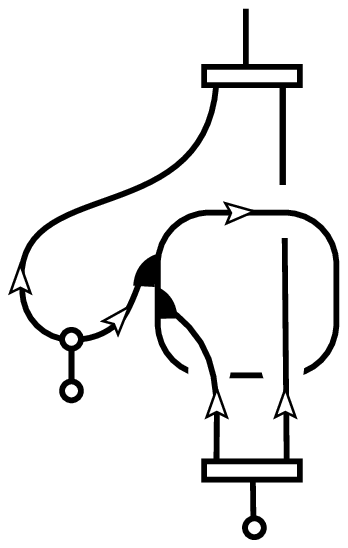}}}
   \put(0,8){
     \setlength{\unitlength}{.75pt}\put(-290,-160){
     \put(337,237)  {\scriptsize $ X{\times} \one $}
     \put(309,273)  {\scriptsize $ A {\times} \one $}
     \put(316,188)  {\scriptsize $ B {\times} \one $}
     \put(376,273)  {\scriptsize $ R $}
     \put(350,320)  {\scriptsize $ Z(A) $}
     \put(330,168)  {\scriptsize $ Z(B) $}
     \put(380,181)  {\scriptsize $ \iota_l $}
     \put(380,294)  {\scriptsize $ r_l $}
  }\setlength{\unitlength}{1pt}}
  \end{picture}}
 = \frac{\dim(X)}{\dim(B)} \,\,\,
  \raisebox{-44pt}{
  \begin{picture}(70,90)
   \put(0,8){\scalebox{.75}{\includegraphics{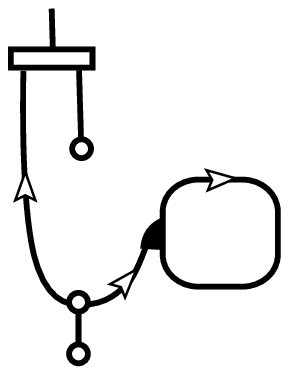}}}
   \put(0,8){
     \setlength{\unitlength}{.75pt}\put(-287,-194){
     \put(313,267)  {\scriptsize $ R $}
     \put(337,257)  {\scriptsize $ X{\times} \one $}
     \put(265,220)  {\scriptsize $ A{\times} \one $}
     \put(291,305)  {\scriptsize $ Z(A) $}
     \put(318,284)  {\scriptsize $ r_l $}
     }\setlength{\unitlength}{1pt}}
  \end{picture}}
 =
  \raisebox{-20pt}{
  \begin{picture}(30,50)
   \put(0,8){\scalebox{.75}{\includegraphics{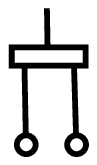}}}
   \put(0,8){
     \setlength{\unitlength}{.75pt}\put(-286,-254){
     \put(290,305)  {\scriptsize $ Z(A) $}
     \put(316,282)  {\scriptsize $ r_l $}
     }\setlength{\unitlength}{1pt}}
  \end{picture}}
\ee
The compatibility with the multiplication,
$m_{Z(A)} \cir (\phi_X \otimes \phi_X) 
= \phi_X \cir m_{Z(B)}$, amounts to the identities
$$
  \left( \frac{\dim(X)}{\dim(B)} \right)^2
 \raisebox{-90pt}{
  \begin{picture}(140,188)
   \put(0,8){\scalebox{.75}{\includegraphics{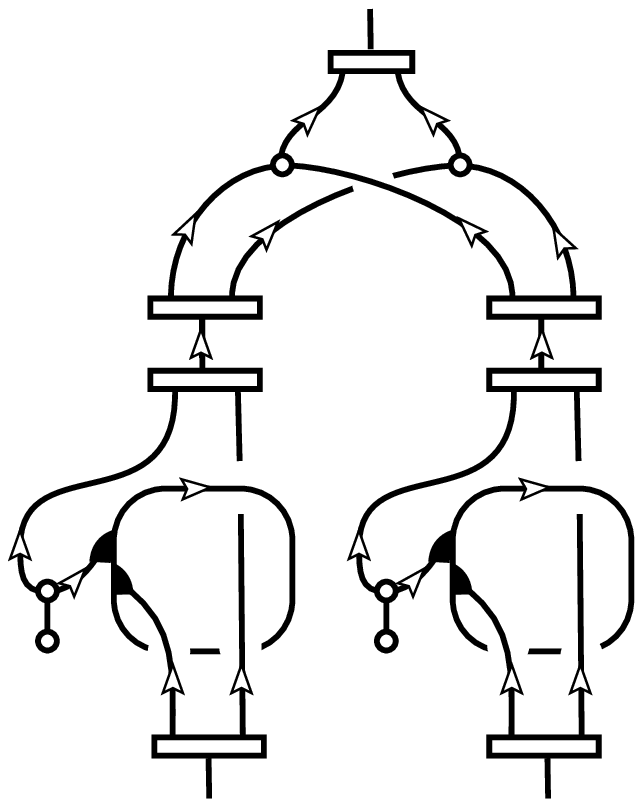}}}
   \put(0,8){
     \setlength{\unitlength}{.75pt}\put(-205,-159){
     \put(252,152)  {\scriptsize $ Z(B)$}
     \put(350,152)  {\scriptsize $ Z(B)$}
     \put(286,173)  {\scriptsize $ \iota_l$}
     \put(382,173)  {\scriptsize $ \iota_l$}
     \put(222,188)  {\scriptsize $ B{\times} \one$}
     \put(279,188)  {\scriptsize $ R $}
     \put(317,188)  {\scriptsize $ B{\times} \one$}
     \put(377,188)  {\scriptsize $ R $}
     \put(239,233)  {\scriptsize $ X{\times} \one$}
     \put(338,233)  {\scriptsize $ X{\times} \one$}
     \put(218,265)  {\scriptsize $ A{\times} \one$}
     \put(315,265)  {\scriptsize $ A{\times} \one$}
     \put(228,290)  {\scriptsize $ Z(A)$}
     \put(325,290)  {\scriptsize $ Z(A)$}
     \put(286,280)  {\scriptsize $ r_l$}
     \put(382,280)  {\scriptsize $ r_l$}
     \put(286,300)  {\scriptsize $ \iota_l$}
     \put(382,300)  {\scriptsize $ \iota_l$}
     \put(220,310)  {\scriptsize $ A{\times} \one$}
     \put(318,310)  {\scriptsize $ A{\times} \one$}
     \put(280,310)  {\scriptsize $ R $}
     \put(377,310)  {\scriptsize $ R $}
     \put(256,355)  {\scriptsize $ A{\times} \one$}
     \put(339,355)  {\scriptsize $ R $}
     \put(327,372)  {\scriptsize $ r_l$}
     \put(299,394)  {\scriptsize $ Z(A)$}
  }\setlength{\unitlength}{1pt}}
  \end{picture}}
 ~\overset{(1)}{=}~ \left( \frac{\dim(X)}{\dim(B)} \right)^2
  \raisebox{-70pt}{
  \begin{picture}(140,148)
   \put(0,8){\scalebox{.75}{\includegraphics{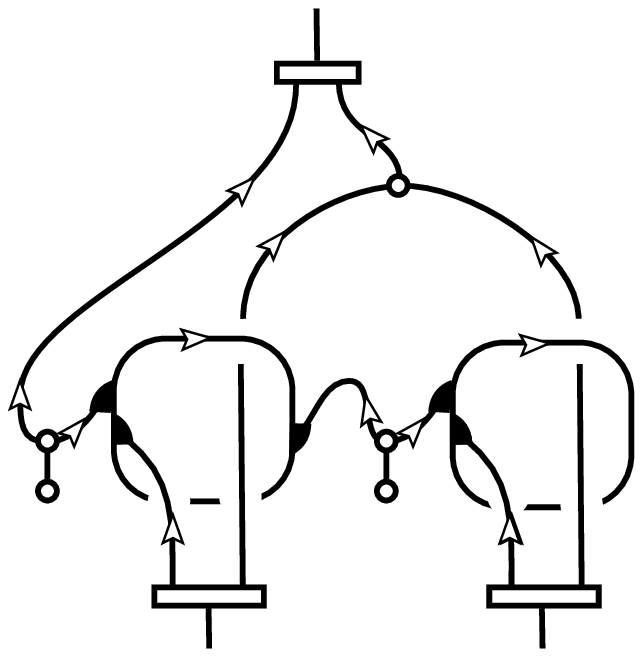}}}
   \put(0,8){
     \setlength{\unitlength}{.75pt}\put(-204,-159){
     \put(252,150)  {\scriptsize $ Z(B)$}
     \put(346,150)  {\scriptsize $ Z(B) $}
     \put(284,174)  {\scriptsize $ \iota_l$}
     \put(380,174)  {\scriptsize $ \iota_l$}
     \put(219,186)  {\scriptsize $ B{\times} \one$}
     \put(316,186)  {\scriptsize $ B{\times} \one$}
     \put(297,245)  {\scriptsize $ A{\times} \one$}
     \put(236,294)  {\scriptsize $ A{\times} \one$}
     \put(278,188)  {\scriptsize $ R $}
     \put(376,188)  {\scriptsize $ R $}
     \put(319,309)  {\scriptsize $ R $}
     \put(239,233)  {\scriptsize $ X{\times} \one$}
     \put(337,233)  {\scriptsize $ X{\times} \one$}
     \put(311,325)  {\scriptsize $ r_l$}
     \put(284,353)  {\scriptsize $ Z(A)$}
    }\setlength{\unitlength}{1pt}}
  \end{picture}}
$$
\be
 \overset{(2)}{=}~
  \frac{\dim(X)^3}{\dim(B)^2\dim(A)}
 \raisebox{-84pt}{
  \begin{picture}(130,176)
   \put(0,8){\scalebox{.75}{\includegraphics{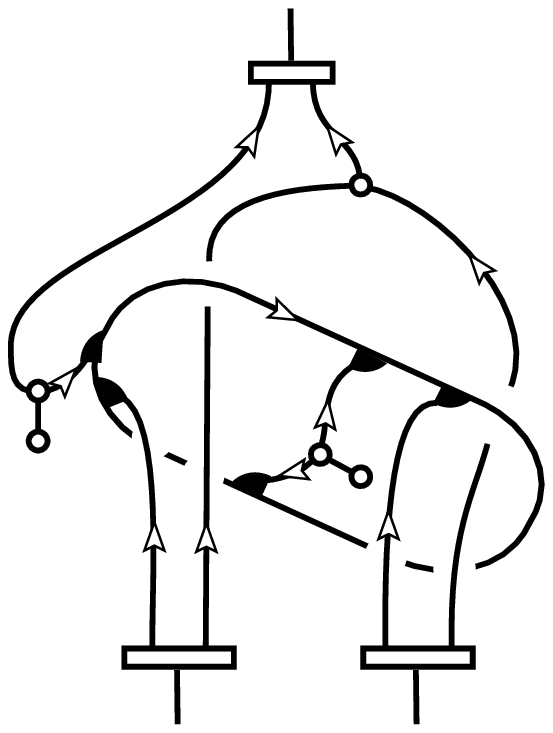}}}
   \put(0,8){
     \setlength{\unitlength}{.75pt}\put(-211,-138){
     \put(251,128)  {\scriptsize $ Z(B) $}
     \put(320,128)  {\scriptsize $ Z(B) $}
     \put(283,157)  {\scriptsize $ \iota_l$}
     \put(352,157)  {\scriptsize $ \iota_l$}
     \put(222,174)  {\scriptsize $ B{\times} \one$}
     \put(289,174)  {\scriptsize $ B{\times} \one$}
     \put(274,222)  {\scriptsize $ B{\times} \one$}
     \put(291,264)  {\scriptsize $ X{\times} \one$}
     \put(275,189)  {\scriptsize $ R $}
     \put(344,169)  {\scriptsize $ R $}
     \put(318,306)  {\scriptsize $ R $}
     \put(247,306)  {\scriptsize $ A{\times} \one$}
     \put(313,324)  {\scriptsize $ r_l$}
     \put(284,352)  {\scriptsize $ Z(A) $}
  }\setlength{\unitlength}{1pt}}
  \end{picture}}
 ~\overset{(3)}{=}~ \frac{\dim(X)}{\dim(B)}
 \raisebox{-84pt}{
\begin{picture}(130,176)
   \put(0,8){\scalebox{.75}{\includegraphics{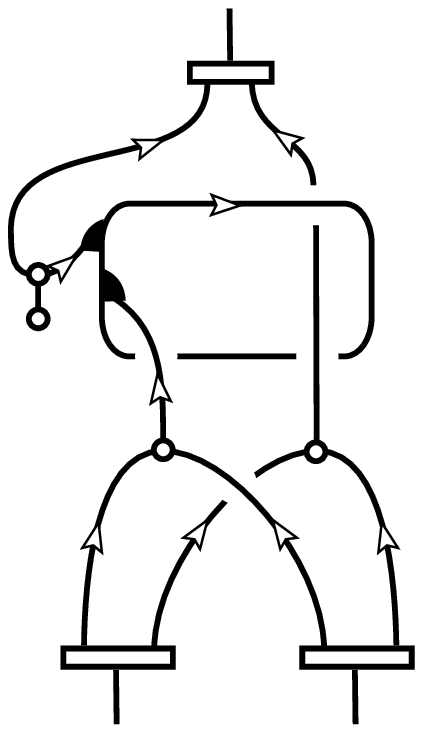}}}
   \put(0,8){
     \setlength{\unitlength}{.75pt}\put(-229,-138){
     \put(251,128)  {\scriptsize $ Z(B) $}
     \put(320,128)  {\scriptsize $ Z(B) $}
     \put(220,170)  {\scriptsize $ B {\times} \one$}
     \put(287,170)  {\scriptsize $ B {\times} \one$}
     \put(241,227)  {\scriptsize $ B {\times} \one$}
     \put(267,191)  {\scriptsize $ R $}
     \put(346,191)  {\scriptsize $ R $}
     \put(316,313)  {\scriptsize $ R $}
     \put(274,273)  {\scriptsize $ X {\times} \one$}
     \put(244,313)  {\scriptsize $ A {\times} \one$}
     \put(312,325)  {\scriptsize $ r_l$}
     \put(284,352)  {\scriptsize $ Z(A) $}
     \put(283,157)  {\scriptsize $ \iota_l$}
     \put(352,157)  {\scriptsize $ \iota_l$}
     }\setlength{\unitlength}{1pt}}
  \end{picture}}
\labl{eq:phiX-alg-map-aux1}
The left hand side is obtained by writing out the definitions of the various morphisms in $m_{Z(A)} \cir (\phi_X \otimes \phi_X)$. In step (1) the two projectors $\iota_l \cir r_l = P_l(R(A))$ have been omitted using lemma \ref{lem:QX-properties}\,(ii,iv), and the uppermost multiplication morphism of $A$ has been replaced by a representation morphism of the bimodule $X$. In step (2) we used property \erf{eq:X-cond-mor}. For step (3) note that the $B \ti \one$-ribbon connecting $X \ti \one$ to itself can be rearranged (using that $B$ is symmetric Frobenius, as well as the representation property) to the projector $P_l(R(B))$ which can be omitted against $\iota_l$. Using the representation property on the remaining two $B\ti\one$-ribbons, as well as \erf{eq:dimX-dimAdimB}, gives the right hand side of \erf{eq:phiX-alg-map-aux1}. Replacing $Q_X = Q_X \cir P_l(R(B)) = Q_X \cir \iota_l \cir r_l$ finally shows that the right hand side is equal to $\phi_X \cir m_{Z(B)}$.
\\[.3em]
e) {\it $\phi_X$ is a coalgebra map}: For this part of the statement, the
coproduct and counit of $Z(A)$ and $Z(B)$ 
have to be normalised as in the proof of
\cite[prop.\,2.37]{corr}. That is, while the multiplication and unit
on $Z(A)$ are given by 
$m_{Z(A)} = r_l \cir m_{R(A)} \cir (\iota_l \oti \iota_l)$ and
$\eta_{Z(A)} = r_l \cir \eta_{R(A)}$, for the coproduct and counit we choose
\be
  \Delta_{Z(A)} = \zeta^{-1}\, 
  (r_l \oti r_l) \cir \Delta_{R(A)} \cir \iota_l
  ~~,~~
  \eps_{Z(A)} = \zeta\,
  \eps_{R(A)} \cir \iota_l
  ~~,\quad
  \zeta = \frac{ \dim(Z(A)) }{\mathrm{Dim}(\Cc)\,\dim(A)}
  ~.
\labl{eq:ZA-coalg}
(That $\dim(Z(A)) \neq 0$ follows from propostion \ref{prop:Z(A)-prop}\,(i).)
In this normalisation one has 
$\eps_{Z(A)} \cir \eta_{Z(A)} = \dim(Z(A))$.
That $\phi_X$ is a coalgebra map can now be verified similarly as in part d) except that at one point one needs the equality between the first and last morphism in the following chain of equalities,
\be
  \raisebox{-42pt}{
  \begin{picture}(70,90)
   \put(0,8){\scalebox{.75}{\includegraphics{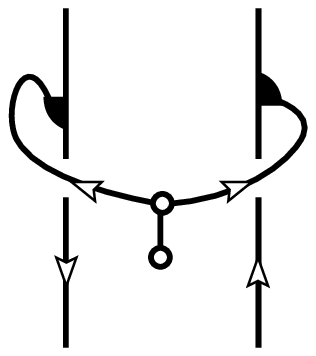}}}
   \put(0,8){
     \setlength{\unitlength}{.75pt}\put(-132,-137){
     \put(143,127)  {\scriptsize $ X^{\vee} $}
     \put(200,127)  {\scriptsize $ X $}
     \put(143,242)  {\scriptsize $ X^{\vee} $}
     \put(200,242)  {\scriptsize $ X $}
     \put(129,185)  {\scriptsize $ B $}
     \put(211,183)  {\scriptsize $ B $}
     }\setlength{\unitlength}{1pt}}
  \end{picture}}
 =
  \raisebox{-52pt}{
  \begin{picture}(130,108)
   \put(0,8){\scalebox{.75}{\includegraphics{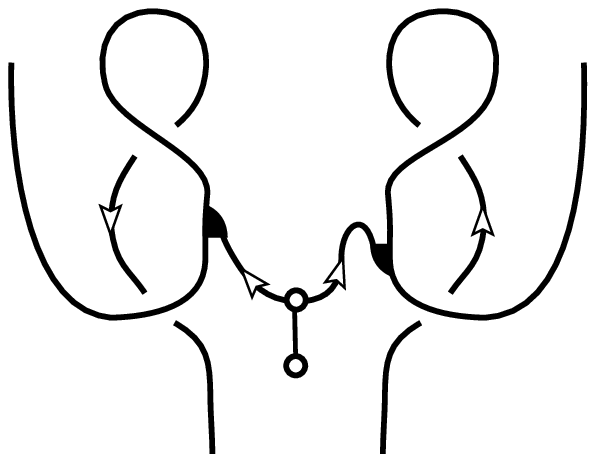}}}
   \put(0,8){
     \setlength{\unitlength}{.75pt}\put(-362,-220){
     \put(415,210)  {\scriptsize $ X^{\vee} $}
     \put(466,210)  {\scriptsize $ X $}
     \put(432,281)  {\scriptsize $ B $}
     \put(445,281)  {\scriptsize $ B $}
     }\setlength{\unitlength}{1pt}}
  \end{picture}}
 = \frac{\dim(A)}{\dim(X)} 
\raisebox{-42pt}{
  \begin{picture}(55,90)
   \put(0,8){\scalebox{.75}{\includegraphics{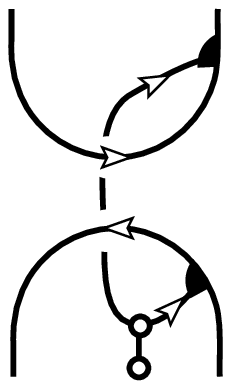}}}
   \put(0,8){
     \setlength{\unitlength}{.75pt}\put(-154,-179){
     \put(152, 169)  {\scriptsize $ X^{\vee} $}
     \put(210, 169)  {\scriptsize $ X $}
     \put(187,271)  {\scriptsize $ A $}
     \put(198,188)  {\scriptsize $ A $}
     \put(152,292)  {\scriptsize $ X^{\vee} $}
     \put(210,292)  {\scriptsize $ X $}
     }\setlength{\unitlength}{1pt}}
  \end{picture}}
  ~~,
\ee
where in the second step \erf{eq:X-cond-mor} is substituted. One also needs to use that $\dim(Z(A)) = \dim(Z(B))$, which follows from part c).
\\[.3em]
This completes the proof of the lemma.
\epf

\begin{lemma}\label{lem:AB-simple-XY}
Let $A$, $B$ be simple non-degenerate algebras and $X$ an 
$A$-$B$-bimodule, $Y$ a $B$-$A$-bimodule such that 
$A \cong X\otimes_B Y$ and $B\cong Y\otimes_A X$ as bimodules. 
Then\\
(i)\phantom{ii} $Y\cong X^{\vee}$ as bimodules,\\
(ii)\phantom{i} $X$ is simple,\\
(iii) the assumptions in lemma \ref{lem:phi_X-iso-frob} hold.
\end{lemma}
\pf
That $Y\cong X^{\vee}$ and that $X$ is simple is proved in 
lemma 3.4 of \cite{defect}. 
Since $A$, $B$, $X$, and $Y$ are all simple as bimodules,
by \cite[lem.\,2.6]{defect} their dimensions are non-zero.
We also have $A \cong X\otimes_B X^\vee$ and $B\cong X^\vee \otimes_A X$ 
as bimodules.
Using this, property \erf{eq:X-cond-mor} follows as a special case
from \cite[eqn.\,(4.8)]{defect}. Property
\erf{eq:X-cond-unit} is proved in lemma 4.1 of \cite{defect}.
\epf

\bigskip\noindent
{\it Proof of (i)$\,\Rightarrow$(ii) in theorem \ref{thm:main}}:
\\[.3em]
By assumption the simple non-degenerate algebras $A$ and $B$ are Morita-equivalent. Therefore there exists an $A$-$B$-bimodule $X$ and a $B$-$A$-bimodule $Y$ such that
$A\cong X\otimes_B Y$ and $B\cong Y \otimes_A X$ as bimodules. Lemma \ref{lem:AB-simple-XY} ensures that the conditions of lemma \ref{lem:phi_X-iso-frob} are met. Thus the morphism $\phi_X: Z(B) \rightarrow Z(A)$ is an isomorphism of algebras. 
\epf

\begin{rema}\label{rem:Z-prop}\end{rema}
\vspace*{-.7em}
(i)~If condition (i) in theorem \ref{thm:main} is met then by proposition \ref{prop:Z(A)-prop}, $Z(A)$ and $Z(B)$ are non-degenerate algebras. The coalgebra structure on $Z(A)$ and $Z(B)$ defined in lemma \ref{lem:nondeg_vs_Frob} is the same as the one used in \erf{eq:ZA-coalg}. Lemma \ref{lem:nondeg_vs_Frob} also implies that $Z(A)$ and $Z(B)$ are even isomorphic as Frobenius algebras.
\\[.3em]
(ii) Recall the definitions of $M_\mathrm{simp}(\Cc)$ and $C_\mathrm{max}(\CxC)$ from remark \ref{rem:intro}\,(ii). The above proof shows that the algebra-isomorphism class $[Z(A)]$ is constant on Morita classes of simple non-degenerate algebras $A$. From proposition \ref{prop:Z(A)-prop}\,(ii) we know that $Z(A)$ is a haploid commutative non-degenerate algebra of dimension $\dim(Z(A)) = \mathrm{Dim}(\Cc)$. Thus $[Z(A)] \In C_\mathrm{max}$. Denoting the Morita class of $A$ by $\{A\}$ it follows that we get a well-defined map $z : M_\mathrm{simp}(\Cc) \rightarrow C_\mathrm{max}(\CxC)$ by setting $z\big( \{A\} \big) = [Z(A)]$, as announced in remark \ref{rem:intro}\,(ii).

\sect{Isomorphic full centre implies Morita equivalence}\label{sec:ii-to-i}

\subsection{The functor $T$} \label{sec:functor-T}

In this section we define a tensor functor $T : \CxC \rightarrow \Cc$ for a braided tensor category $\Cc$. For concreteness, we will spell out associators and unit constraints explicitly.
The monoidal structure on $\Cc$ consists of 
the unit object $\one$ and the tensor-product bifunctor
$\otimes: \Cc\ti \Cc \rightarrow \Cc$, 
together with a left unit isomorphism 
$l_U: \one \oti U \rightarrow U$,
a right unit isomorphism $r_U: U\oti \one \rightarrow U$ for
each $U \In \Cc$, and an associator
$\alpha_{U,V,W}: U\otimes (V\otimes W) 
\rightarrow (U\otimes V) \otimes W$ for any triple objects 
$U,V,W\in \Cc$. 

The bifunctor $\otimes$ can be naturally extended to a functor 
$T: \CxC \rightarrow \Cc$. 
Namely, 
$T( \oplus_{i=1}^{N} U_i \times V_i) = \oplus_{i=1}^N U_i \otimes V_i$
for all $U_i, V_i \in \Cc$ and $N\in \Nb$. 
Let $\varphi_0: \one \rightarrow T(\one \ti \one)$ be 
$l_{\one}^{-1}$. For $U, V, W, X\in \Cc$, notice that
\be\bearll
T\big(U\ti V\big) \otimes T\big(W\ti X\big) 
\etb= (U\oti V) \otimes (W\oti X), \enl 
T\big( (U\ti V) \oti (W\ti X) \big) \etb= 
(U\oti W) \otimes (V\oti X).
\eear\ee
We define $\varphi_2: T\big(U\times V\big) \otimes T\big(W\times X\big)
\rightarrow T\big( (U\times V) \otimes (W\times X)\big)$ by
\be
\varphi_2 := \alpha_{U,W,V\otimes X}
\circ \big(\id_U \otimes \alpha_{W,V,X}^{-1}\big)
\circ \big(\id_U \otimes (c_{WV}^{-1} \otimes \id_X)\big) \circ 
(\id_U \otimes \alpha_{V,W,X}) \circ \alpha_{U,V,W\otimes X}^{-1}.
\ee
The above definition of $\varphi_2$ can be naturally extended to a morphism 
$T(M_1)\otimes T(M_2) \rightarrow T(M_1 \otimes M_2)$ for 
any pair of objects $M_1, M_2$ in $\CxC$. 
We still denote the extended morphism as $\varphi_2$. 
(We hide the dependence of $\varphi_2$ 
on $M_1, M_2$ in our notation for simplicity.)
We have

\begin{lemma}\label{lem:T-monoidal}
The functor $T$ together with 
$\varphi_0$ and $\varphi_2$ 
is a tensor functor. 
\end{lemma}

Note that $T$ takes algebras to algebras (see for example 
\cite[prop.\,3.7]{ko-ocfa}) 
but in general does not preserve commutativity. 
Explicitly, if $(B, m_{B}, \eta_B)$ is an algebra in 
$\mathcal{C} \boxtimes \tilde{\mathcal{C}}$, then the triple
$(T(B), m_{T(B)}, \eta_{T(B)})$, where 
\be 
m_{T(B)}: = T(m_B) \circ \varphi_2~~, \qquad 
\eta_{T(B)} := T(\eta_B) \circ \varphi_0~,
\labl{eq:TB}
is an algebra in $\mathcal{C}$.

\subsection{The full centre transported to $\Cc$ and simple modules}
\label{sec:fc-simpmod}

Let now $\Cc$ again be a (strict) modular tensor category, and let $A$ be a simple non-degenerate algebra in $\Cc$. As observed in section \ref{sec:intro}, this implies in particular that $\dim(A) \neq 0$. The category of left $A$-modules is again semisimple and abelian \cite[props.\,5.1 and 5.24]{fs-cat} with a finite number of isomorphism classes of simple objects (this follows e.g.\ by combining the fact that $\Cc$ itself only has a finite number of isomorphism classes of simple objects with \cite[lem.\,4.15]{fs-cat}). Let $\{\, M_\kappa \,|\, \kappa \In \Jc \,\}$ be a set of representatives of the isomorphism classes of simple left $A$-modules.

\begin{lemma}\label{lem:MxM-prop1}
Let $A$ be a non-degenerate algebra in $\Cc$ and let $M$ be a left $A$-module. 
\\[.3em]
(i)\phantom{i} $M^\vee \otimes_A M$ is an algebra with unit $e_A \cir \tilde b_M$ and multiplication $r_A \cir (\id_{M^\vee} \oti \tilde d_M \oti \id_M) \circ (e_A \otimes e_A)$.
\\[.3em]
(ii) $M$ is simple if and only if $M^\vee \otimes_A M$ is haploid.
\end{lemma}
\pf Part (i) is a straightforward calculation, see e.g.\ \cite[eqn.\,(2.48)]{tft2}. Claim (ii) follows since $\Hom_A(M,M) \cong \Hom(M^\vee \otimes_A M,\one)$. The first space is one-dimensional iff $M$ is simple, and the second space is one-dimensional iff $M^\vee \otimes_A M$ is haploid.
\epf

\medskip

We define two algebras $C_A$ and $T_A$ in $\Cc$ as follows,
\be
  C_A ~=~ T\big(Z(A)\big) \quad , \qquad
  T_A ~=~ \bigoplus_{\kappa \in \Jc} M_\kappa^\vee \otimes_A M_\kappa
  ~.
\labl{eq:CA-TA-def}
{}From the discussion in section \ref{sec:functor-T} we see that 
$C_A$ is naturally an algebra in $\Cc$, and by lemma \ref{lem:MxM-prop1}
the same holds for $T_A$. Note that $C_A$ is not necessarily 
commutative, even though $Z(A)$ is. 

\begin{prop} \label{prop:CA=TA}
$C_A \cong T_A$ as algebras. 
\end{prop}

As an isomorphism between objects, rather than algebras, this statement can already be found in the conformal field theory literature, see \cite[eqn.\,(4.2)]{ss}.

The proof of proposition \ref{prop:CA=TA} needs a bit of preparation and will be given at the end of this section. We start by recalling the definition of local morphisms in $\Hom(A \oti U,V)$ from \cite[sect.\,5.3]{tft1}. Define the morphism $P^l_A(U) : A \oti U \rightarrow A \oti U$ as
\be
  P^l_A(U) ~=~ 
  \raisebox{-40pt}{
  \begin{picture}(65,80)
   \put(0,8){\scalebox{.75}{\includegraphics{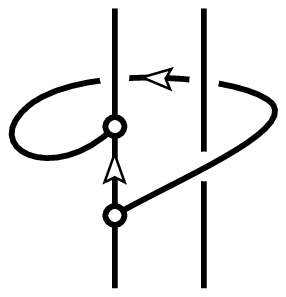}}}
   \put(0,8){
     \setlength{\unitlength}{.75pt}\put(-18,-19){
     \put(45, 10)  {\scriptsize $ A $}
     \put(45,105)  {\scriptsize $ A $}
     \put(70, 10)  {\scriptsize $ U $}
     \put(70,105)  {\scriptsize $ U $}
     \put(93, 79)  {\scriptsize $ A $}
     \put(55, 65)  {\scriptsize $ m $}
     \put(36, 36)  {\scriptsize $ \Delta $}
     }\setlength{\unitlength}{1pt}}
  \end{picture}}
  ~~.
\ee
One verifies that $P^l_A(U)$ is an idempotent, cf.\ \cite[lem.\,5.2]{tft1}. Note that the idempotent defining the left centre can be written as $P_l(A) = P^l_A(\one)$. We set
\bea
  \Hom_\mathrm{loc}(A \oti U,V)
  ~=~ 
  \big\{ \, f : A \oti U \rightarrow V \,\big|\, 
  f \cir P^l_A(U) = f \, \big\} ~.
\eear\labl{eq:homloc-def}
The morphisms in $\Hom_\mathrm{loc}(A \oti U,V)$ are called 
{\em local}. Let $\{ \, \mu^i_\alpha \, \}$ be a basis of 
$\Hom(A \oti U_i,U_i)$ such that $\mu^i_\alpha$ is local 
for $\alpha=1,\dots,N_i^\mathrm{loc}$ and $\mu^i_\alpha \cir P^l_A(U_i) = 0$ 
for $\alpha > N_i^\mathrm{loc}$. Let $\{ \, \bar\mu^i_\alpha \, \}$ be the 
basis of $\Hom(U_i,A \oti U_i)$ that is dual to $\mu^i_\alpha$ 
in the sense that $\mu^i_\alpha
\cir \bar\mu^i_\beta = \delta_{\alpha,\beta}\, \id_{U_i}$.

One can prove that $N_i^\mathrm{loc} = \dim\Hom(Z(A),U_i \ti U_i^\vee)$,
see \cite[lem.\,5.6]{tft1}. By propostion \ref{prop:Z(A)-prop}\,(ii), 
$Z(A)$ is haploid and so $N_0^\mathrm{loc}=1$. Let us agree to
choose the basis vector in $\Hom_\mathrm{loc}(A,\one)$ to be
$\mu^0_1 = \dim(A)^{-1} \eps_A$, and consequently also 
$\bar\mu^0_1 = \eta_A$.

Using these bases of local morphisms we can define numbers 
$s^A_{\kappa,i\alpha}$ and $\tilde s^A_{i\alpha,\kappa}$ as
in \cite[sect.\,5.7]{tft1},
\be
  s^A_{\kappa,i\alpha} ~=~
  \raisebox{-35pt}{
  \begin{picture}(100,78)
   \put(0,0){\scalebox{.75}{\includegraphics{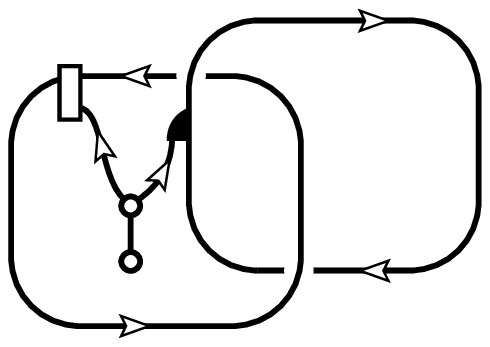}}}
   \put(0,0){
     \setlength{\unitlength}{.75pt}\put(-21,-19){
     \put(40,60)  {\scriptsize $ A $}
     \put(76,63)  {\scriptsize $ M_\kappa $}
     \put(34,104)  {\scriptsize $ \mu^i_\alpha $}
     \put(109,69)  {\scriptsize $ U_i $}
     }\setlength{\unitlength}{1pt}}
  \end{picture}}
  \quad , \qquad
  \tilde s^A_{i\alpha,\kappa} ~=~
  \raisebox{-35pt}{
  \begin{picture}(100,78)
   \put(0,0){\scalebox{.75}{\includegraphics{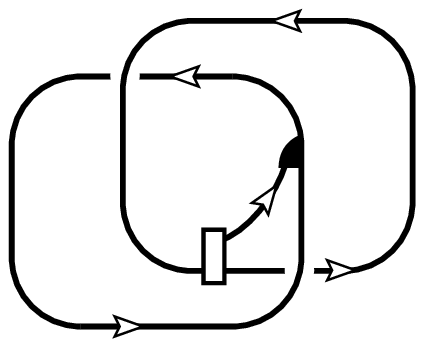}}}
   \put(0,0){
     \setlength{\unitlength}{.75pt}\put(-22,-21){
     \put(110,62)  {\scriptsize $ M_\kappa $}
     \put(74,59)  {\scriptsize $ \bar\mu^i_\alpha $}
     \put(142,76)  {\scriptsize $ U_i $}
     }\setlength{\unitlength}{1pt}}
  \end{picture}}
  ~~,
\ee
where $\kappa \In \Jc$, $i \in \Ic$ and 
$\alpha=1,\dots,N^\mathrm{loc}_i$. We have
\be
\sum_{i\in\Ic} \sum_{\alpha=1}^{N^\mathrm{loc}_i}
s^A_{\kappa,i\alpha} \, \tilde s^A_{i\alpha,\lambda}
= \mathrm{Dim}(\Cc) \, \delta_{\kappa,\lambda}
\quad , \quad
\sum_{\kappa \in \Jc}
\tilde s^A_{i\alpha,\kappa} \, s^A_{\kappa,j\beta} 
= \mathrm{Dim}(\Cc) \, \delta_{i,j} \, \delta_{\alpha,\beta}
~.
\labl{eq:sA-stildeA}
where in the first equality $\kappa,\lambda \in \Jc$ and in the 
second equality $i,j \in \Ic$ have to be chosen such that
$N^\mathrm{loc}_i >0$ and $N^\mathrm{loc}_j >0$.
These equalities are proved in
\cite[prop.\,5.16 and 5.17]{tft1}. 
They imply in particular that $s^A$ and $\tilde s^A$ are
square matrices, $\sum_{i \in \Ic} N^\mathrm{loc}_i = |\Jc|$.
We are now in a position to prove the following lemma.

\begin{lemma}\label{lem:M-loop}
\be
  \sum_{\kappa \in \Jc} \frac{\dim(M_\kappa)}{\mathrm{Dim}(\Cc)}
  \raisebox{-50pt}{
  \begin{picture}(65,103)
   \put(0,8){\scalebox{.75}{\includegraphics{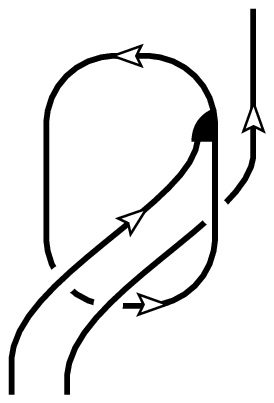}}}
   \put(0,8){
     \setlength{\unitlength}{.75pt}\put(-20,-13){
     \put(80,46)  {\scriptsize $ M_\kappa $}
     \put(17,3)  {\scriptsize $ A $}
     \put(33,3)  {\scriptsize $ U_i $}
     \put(87,130)  {\scriptsize $ U_i $}
     }\setlength{\unitlength}{1pt}}
  \end{picture}} 
  ~=~ \delta_{i,0} \, \eps_A ~~.
\labl{eq:Mloop}
\end{lemma}
\pf
Denote the left hand side of \erf{eq:Mloop} by $f$ and the morphism represented pictorially by $f_\kappa$, 
s.t.\ $f = \mathrm{Dim}(\Cc)^{-1} \sum_{\kappa} \dim(M_\kappa) f_\kappa$.
A calculation similar to the one needed to show that $P^l_A(U_i)$ is an idempotent shows that $f_\kappa \cir P^l_A(U_i) = f_\kappa$. 
Thus also $f \cir P^l_A(U_i) = f$ and hence 
$f \in \Hom_\mathrm{loc}(A \oti U_i,U_i)$. We can therefore
expand $f$ in the basis $\mu^i_\alpha$ as
$f = \sum_{\beta=1}^{N^\mathrm{loc}_i} c_\beta \,\mu^i_\beta$. To
determine the constants $c_\beta$ we compose both sides with the
dual basis element $\bar\mu^i_\alpha$ from the right. This results
in 
$\mathrm{Dim}(\Cc)^{-1}
\sum_{\kappa} \dim(M_\kappa)
f_\kappa \cir\bar \mu^i_\alpha
= c_\alpha \id_{U_i}$.
The constant $c_\alpha$ can then be extracted by taking the trace
on both sides, 
\be\bearll
  c_\alpha \etb=  
  \frac{1}{\dim(U_i)\mathrm{Dim}(\Cc)} 
  \sum_{\kappa \in \Jc}
  \dim(M_\kappa) \,
  \mathrm{tr}_{U_i}\big( f_\kappa \cir\bar \mu^i_\alpha \big)
  \enl
  \etb
  =
  \frac{\dim(A)}{\dim(U_i)\mathrm{Dim}(\Cc)} 
  \sum_{\kappa \in \Jc} s^A_{\kappa,01} \,
  \tilde s^A_{i\alpha,\kappa} 
  = \dim(A) \, \delta_{i,0} \, \delta_{\alpha,1}~,
\eear\ee
where in the second step we used that
$s^A_{\kappa,01} = \dim(M_\kappa)/\dim(A)$ (recall the choice $\mu^0_1 = \dim(A)^{-1} \eps_A$) and 
$\mathrm{tr}_{U_i}\big( f_\kappa \cir\bar \mu^i_\alpha \big)
= \tilde s^A_{i\alpha,\kappa}$ which follows by comparing
the pictorial representations of the morphisms on either side.
The third step is a consequence of the second equality in
\erf{eq:sA-stildeA}. Substituting this result for $c_\alpha$
back into 
$f = \sum_{\beta=1}^{N^\mathrm{loc}_i} c_\beta \, \mu^i_\beta$
then yields \erf{eq:Mloop}.
\epf

\medskip

We will also need the following identity.

\begin{lemma}\label{lem:MM-Ui-loop}
\be  
\sum_{i \in \Ic} \dim(U_i)
  \raisebox{-48pt}{
  \begin{picture}(60,100)
   \put(0,8){\scalebox{.75}{\includegraphics{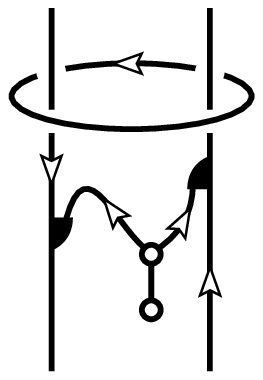}}}
   \put(0,8){
     \setlength{\unitlength}{.75pt}\put(-15,-14){
     \put(23,  5)  {\scriptsize $ M_\alpha^\vee $}
     \put(23,125)  {\scriptsize $ M_\alpha^\vee $}
     \put(68,  5)  {\scriptsize $ M_\beta $}
     \put(68,125)  {\scriptsize $ M_\beta $}
     \put(49,66)  {\scriptsize $ A $}
     \put(47,111)  {\scriptsize $ U_i $}
     }\setlength{\unitlength}{1pt}}
  \end{picture}} 
~=~
\delta_{\alpha,\beta} \, \frac{\mathrm{Dim}(\Cc)}{\dim(M_\alpha)}
  \raisebox{-48pt}{
  \begin{picture}(60,100)
   \put(0,8){\scalebox{.75}{\includegraphics{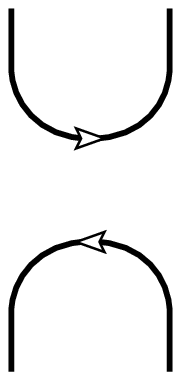}}}
   \put(0,8){
     \setlength{\unitlength}{.75pt}\put(-26,-14){
     \put(23,  5)  {\scriptsize $ M_\alpha^\vee $}
     \put(23,125)  {\scriptsize $ M_\alpha^\vee $}
     \put(68,  5)  {\scriptsize $ M_\alpha $}
     \put(68,125)  {\scriptsize $ M_\alpha $}
     }\setlength{\unitlength}{1pt}}
  \end{picture}}
  . 
\labl{eq:MM-Ui-loop}
\end{lemma}
\pf
Let $\{\,x_\nu^k\,\}$ be a basis of 
$\Hom_A(M_{\alpha}\otimes U_k,  M_{\beta})$
and let $\{\,\bar x_\nu^k\,\}$ be the basis 
of $\Hom_A(M_{\beta}, M_{\alpha}\otimes U_k)$ dual
to $x_\nu^k$ in the sense that 
$ x_\mu^k \cir \bar x_\nu^k = \delta_{\mu,\nu} \id_{M_\beta}$.
For $k=0$ and $\alpha=\beta$
there is only one basis vector in each space,
and we choose $x_1^0 = \bar x_1^0 = \id_{M_\alpha}$.

Using the identity \cite[eqn\,(4.8)]{defect} 
(actually we need the `vertically reflected' version) in the
special case of $A$-$\one$-bimodules, we obtain 
\be 
\sum_{i \in \Ic}  \dim(U_i)
  \raisebox{-48pt}{
  \begin{picture}(60,100)
   \put(0,8){\scalebox{.75}{\includegraphics{pic-MM-Ui-loop1.eps}}}
   \put(0,8){
     \setlength{\unitlength}{.75pt}\put(-15,-14){
     \put(23,  5)  {\scriptsize $ M_\alpha^\vee $}
     \put(23,125)  {\scriptsize $ M_\alpha^\vee $}
     \put(68,  5)  {\scriptsize $ M_\beta $}
     \put(68,125)  {\scriptsize $ M_\beta $}
     \put(49,66)  {\scriptsize $ A $}
     \put(47,111)  {\scriptsize $ U_i $}
     }\setlength{\unitlength}{1pt}}
  \end{picture}} 
~=~
\sum_{i,k,\nu} \dim(U_i) \,  \frac{\dim(U_k)}{\dim(M_\beta)}
  \raisebox{-48pt}{
  \begin{picture}(70,100)
   \put(0,8){\scalebox{.75}{\includegraphics{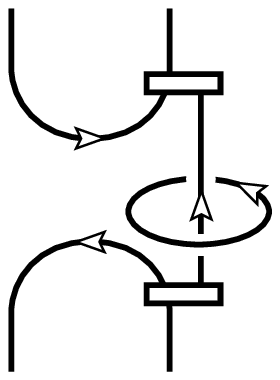}}}
   \put(0,8){
     \setlength{\unitlength}{.75pt}\put(-26,-14){
     \put(23,  5)  {\scriptsize $ M_\alpha^\vee $}
     \put(23,125)  {\scriptsize $ M_\alpha^\vee $}
     \put(68,  5)  {\scriptsize $ M_\beta $}
     \put(68,125)  {\scriptsize $ M_\beta $}
     \put(92,98)  {\scriptsize $ x_\nu^k $}
     \put(92,36)  {\scriptsize $ \bar x_\nu^k $}
     \put(104,59)  {\scriptsize $ U_i $}
     \put(85,83)  {\scriptsize $ U_k $}
     }\setlength{\unitlength}{1pt}}
  \end{picture}} ~~.
\labl{eq:MM-Ui-loop-aux1}
Using further \erf{eq:Mloop} in the special case $A=\one$
(or directly eqn.\,(3.1.19) in \cite{baki}) one finds that
the right hand side of \erf{eq:MM-Ui-loop-aux1} is equal to
\be
  \sum_{k} \frac{\dim(U_k)}{\dim(M_\beta)}  \,
  \mathrm{Dim}(\Cc) \, \delta_{k,0} \, \delta_{\alpha,\beta} ~
  \tilde b_{M_\alpha} \cir d_{M_\alpha} ~,
\ee
which in turn is equal to the right hand side of \erf{eq:MM-Ui-loop}.
\epf

\medskip

For $\iota_l : Z(A) \rightarrow R(A)$ and $r_l : R(A) \rightarrow Z(A)$ the embedding and restriction morphisms of the full centre as in section \ref{sec:pre-full}, let
\be
  e_C = T(\iota_l) : C_A \rightarrow T(R(A)) \qquad \text{and} \qquad
  r_C = T(r_l) : T(R(A)) \rightarrow C_A ~.
\labl{eq:ec-rc-def}
Note that $T(R(A)) = \bigoplus_{i \in \Ic} A \oti U_i \oti U_i^\vee$. Let further 
$e_i : C_A \rightarrow A \oti U_i \oti U_i^\vee$ 
and
$r_i : A \oti U_i \oti U_i^\vee \rightarrow C_A$ 
be given by the compositions
\be
e_i = C_A \overset{e_C}\hookrightarrow T(R(A)) 
\twoheadrightarrow A \oti U_i \oti U_i^\vee
\quad \text{and} \quad
r_i = A \oti U_i \oti U_i^\vee \hookrightarrow T(R(A)) 
\overset{r_C}{\twoheadrightarrow} C_A
~.
\labl{eq:ei-ri-def_CA}
For $T_A = \bigoplus_{\kappa\in\Jc} M_\kappa^\vee \otimes_A M_\kappa$ 
we define in the same way
$e_\kappa : T_A \rightarrow M_\kappa^\vee \oti M_\kappa$ 
and
$r_\kappa : M_\kappa^\vee \oti M_\kappa \rightarrow T_A$ 
to be the compositions
\be
e_\kappa = T_A \twoheadrightarrow 
M_\kappa^\vee \otimes_A M_\kappa
\overset{e_A}{\hookrightarrow} M_\kappa^\vee \oti M_\kappa
\quad \text{and} \quad
r_\kappa = M_\kappa^\vee \oti M_\kappa 
\overset{r_A}{\twoheadrightarrow} 
M_\kappa^\vee \otimes_A M_\kappa \hookrightarrow T_A
~.
\labl{eq:e_kap-r_kap}
Using these ingredients we
define two morphisms $\varphi: C_A \rightarrow T_A$ 
and $\bar{\varphi}: T_A \rightarrow C_A$ by 
\be
\varphi ~=~ \sum_{i \in \Ic} \sum_{\kappa \in \Jc}
  \raisebox{-58pt}{
  \begin{picture}(64,120)
   \put(0,8){\scalebox{.75}{\includegraphics{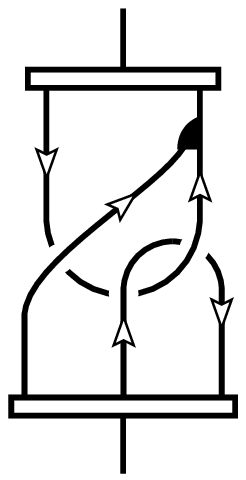}}}
   \put(0,8){
     \setlength{\unitlength}{.75pt}\put(-18,-0){
     \put(77,86)  {\scriptsize $ M_\kappa $}
     \put(38,32)  {\scriptsize $ U_i $}
     \put(45,85)  {\scriptsize $ A $}
     \put(83,115)  {\scriptsize $ r_\kappa $}
     \put(87,19)  {\scriptsize $ e_i $}
     \put(45,140)  {\scriptsize $ T_A $}
     \put(45,-9)  {\scriptsize $ C_A $}
     }\setlength{\unitlength}{1pt}}
  \end{picture}} 
\quad \text{and} \quad
\bar\varphi ~=~ \sum_{i \in \Ic} \sum_{\kappa \in \Jc}
\frac{\dim(U_i)\dim(M_\kappa)}{\mathrm{Dim}(\Cc)}
  \raisebox{-58pt}{
  \begin{picture}(64,120)
   \put(0,8){\scalebox{.75}{\includegraphics{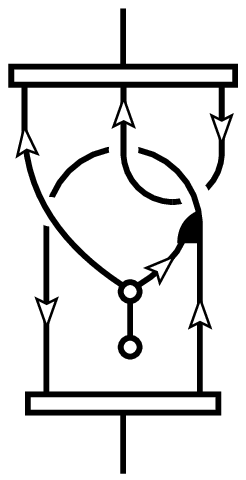}}}
   \put(0,8){
     \setlength{\unitlength}{.75pt}\put(-18,-0){
     \put(77,56)  {\scriptsize $ M_\kappa $}
     \put(56,101)  {\scriptsize $ U_i $}
     \put(47,65)  {\scriptsize $ A $}
     \put(87,115)  {\scriptsize $ r_i $}
     \put(83,19)  {\scriptsize $ e_\kappa $}
     \put(45,140)  {\scriptsize $ C_A $}
     \put(45,-9)  {\scriptsize $ T_A $}
     }\setlength{\unitlength}{1pt}}
  \end{picture}} 
  .
\labl{eq:phi-phi-def}

\begin{lemma}\label{lem:phi-inv-1}
$\varphi \circ \bar{\varphi} = \id_{T_A}$. 
\end{lemma}
\pf
Let $c_{i\lambda}:= \dim(U_i)\dim(M_{\lambda}) / \mathrm{Dim}(\Cc)$. 
Consider the equalities
$$
\varphi \circ \bar{\varphi}
\overset{(1)}{=} \sum_{i,\kappa,\lambda} c_{i\lambda} 
  \raisebox{-105pt}{
  \begin{picture}(74,210)
   \put(0,8){\scalebox{.75}{\includegraphics{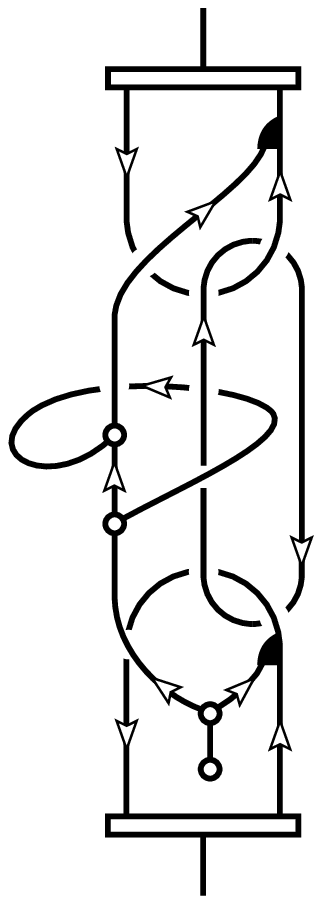}}}
   \put(0,8){
     \setlength{\unitlength}{.75pt}\put(-13,-10){
     \put(64,0)  {\scriptsize $ T_A $}
     \put(64,272)  {\scriptsize $ T_A $}
     \put(101,30)  {\scriptsize $ e_\lambda $}
     \put(101,246)  {\scriptsize $ r_\kappa $}
     \put(95,220)  {\scriptsize $ M_\kappa $}
     \put(95,63)  {\scriptsize $ M_\lambda $}
     \put(71,114)  {\scriptsize $ U_i $}
     \put(33,98)  {\scriptsize $ A $}
     \put(33,163)  {\scriptsize $ A $}
     }\setlength{\unitlength}{1pt}}
  \end{picture}} 
\overset{(2)}{=} \sum_{i,\kappa,\lambda} c_{i\lambda} 
  \raisebox{-105pt}{
  \begin{picture}(74,210)
   \put(0,8){\scalebox{.75}{\includegraphics{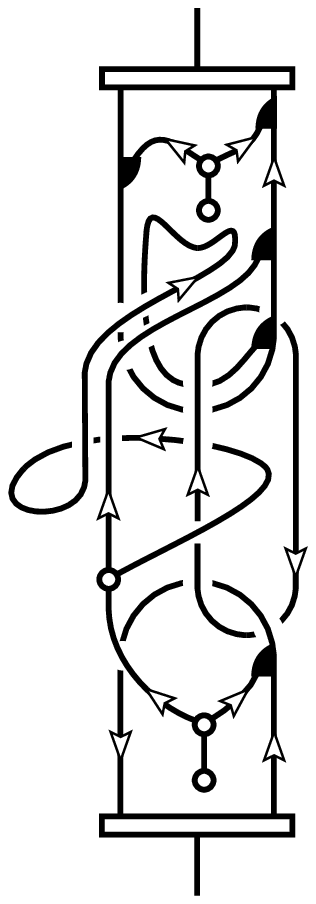}}}
   \put(0,8){
     \setlength{\unitlength}{.75pt}\put(-13,-10){
     \put(64,0)  {\scriptsize $ T_A $}
     \put(64,272)  {\scriptsize $ T_A $}
     \put(101,30)  {\scriptsize $ e_\lambda $}
     \put(101,246)  {\scriptsize $ r_\kappa $}
     \put(95,220)  {\scriptsize $ M_\kappa $}
     \put(95,63)  {\scriptsize $ M_\lambda $}
     \put(55,122)  {\scriptsize $ U_i $}
     \put(32,109)  {\scriptsize $ A $}
     \put(25,163)  {\scriptsize $ A $}
     \put(64,232)  {\scriptsize $ A $}
     }\setlength{\unitlength}{1pt}}
  \end{picture}} 
\overset{(3)}{=} \sum_{i,\kappa,\lambda} c_{i\lambda} 
  \raisebox{-105pt}{
  \begin{picture}(60,210)
   \put(0,8){\scalebox{.75}{\includegraphics{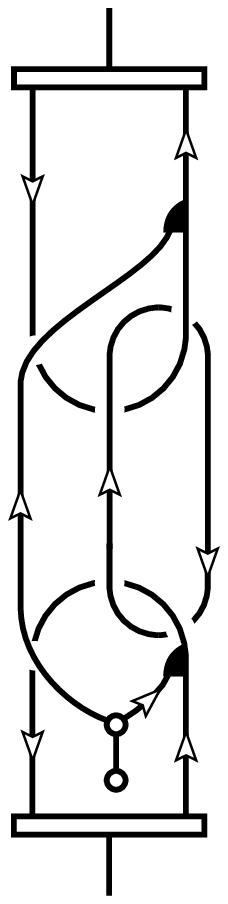}}}
   \put(0,8){
     \setlength{\unitlength}{.75pt}\put(-33,-10){
     \put(64,0)  {\scriptsize $ T_A $}
     \put(64,272)  {\scriptsize $ T_A $}
     \put(101,30)  {\scriptsize $ e_\lambda $}
     \put(101,246)  {\scriptsize $ r_\kappa $}
     \put(95,220)  {\scriptsize $ M_\kappa $}
     \put(95,63)  {\scriptsize $ M_\lambda $}
     \put(55,122)  {\scriptsize $ U_i $}
     \put(32,109)  {\scriptsize $ A $}
     }\setlength{\unitlength}{1pt}}
  \end{picture}} 
$$
\be
\overset{(4)}{=} \sum_{i,\kappa,\lambda} c_{i\lambda} 
  \raisebox{-85pt}{
  \begin{picture}(115,177)
   \put(0,8){\scalebox{.75}{\includegraphics{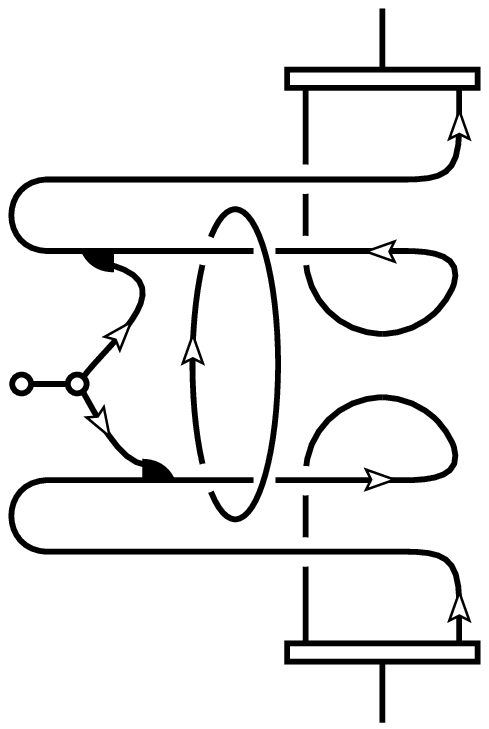}}}
   \put(0,8){
     \setlength{\unitlength}{.75pt}\put(-7,-0){
     \put(110,-8)  {\scriptsize $ T_A $}
     \put(110,214)  {\scriptsize $ T_A $}
     \put(147,20)  {\scriptsize $ e_\lambda $}
     \put(147,186)  {\scriptsize $ r_\kappa $}
     \put(140,126)  {\scriptsize $ M_\kappa $}
     \put(140,78)  {\scriptsize $ M_\lambda $}
     \put(66,103)  {\scriptsize $ U_i $}
     \put(26,113)  {\scriptsize $ A $}
     }\setlength{\unitlength}{1pt}}
  \end{picture}} 
\overset{(5)}{=} \sum_{\kappa,\lambda} \delta_{\kappa,\lambda}
  \raisebox{-85pt}{
  \begin{picture}(100,177)
   \put(0,8){\scalebox{.75}{\includegraphics{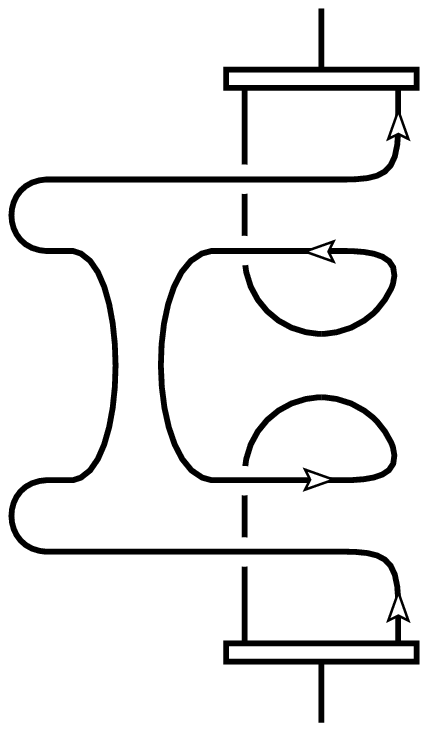}}}
   \put(0,8){
     \setlength{\unitlength}{.75pt}\put(-24,-0){
     \put(110,-8)  {\scriptsize $ T_A $}
     \put(110,214)  {\scriptsize $ T_A $}
     \put(147,20)  {\scriptsize $ e_\kappa $}
     \put(147,186)  {\scriptsize $ r_\kappa $}
     \put(36,103)  {\scriptsize $ M_\kappa $}
     \put(140,78)  {\scriptsize $ M_\kappa $}
     }\setlength{\unitlength}{1pt}}
  \end{picture}} 
~~.
\labl{eq:phi-inv-1-aux}
Step (1) amounts to the definition of $\varphi$ and $\bar\varphi$
and to the identity $e_i \cir r_i = P^l_A(U_i) \oti \id_{U_i^\vee}$.
Step (2) and (3) show that the idempotent $P^l_A(U_i)$ can be 
cancelled against $r_\kappa$. To this end $r_\kappa$ is replaced
by $r_\kappa \cir P_{\otimes A}$ and the multiplication morphism
is moved to the $M_\kappa$-ribbon, as indicated. In doing so one
uses that $A$ is symmetric Frobenius and that $M_\kappa$ is a
left $A$-module. In step (3) one uses the representation property
once more, as well as the fact that $A$ is normalised-special.
Step (4) is just a deformation of the ribbon graph so that one
can apply lemma \ref{lem:MM-Ui-loop}. This is done in step (5), and after
`straightening' the $M_\kappa$-ribbons and using 
$\sum_{\kappa \in \Jc} r_\kappa \cir e_\kappa = \id_{T_A}$, one
finally obtains that the right hand side of \erf{eq:phi-inv-1-aux}
is equal to $\id_{T_A}$.
\epf

\begin{lemma}\label{lem:phi-inv-2}
$\bar{\varphi} \circ \varphi = \id_{C_A}$
\end{lemma}
\pf As in the proof of the previous lemma we set
$c_{i\kappa}:= \dim(U_i)\dim(M_{\kappa}) / \mathrm{Dim}(\Cc)$.
Consider the equalities
$$
\bar \varphi \circ \varphi 
\overset{(1)}{=} \sum_{i,j,\kappa} c_{i\kappa} 
  \raisebox{-91pt}{
  \begin{picture}(60,190)
   \put(0,8){\scalebox{.75}{\includegraphics{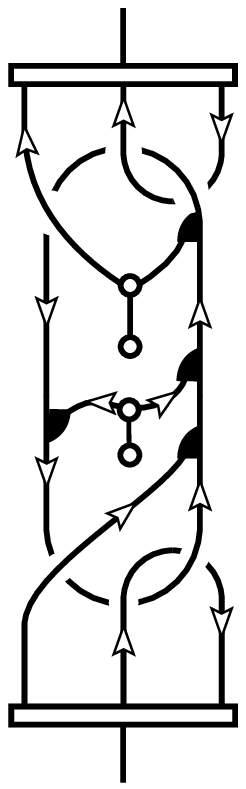}}}
   \put(0,8){
     \setlength{\unitlength}{.75pt}\put(-17,0){
     \put(47,-10)  {\scriptsize $ C_A $}
     \put(47,230)  {\scriptsize $ C_A $}
     \put(10,173)  {\scriptsize $ A $}
     \put(56,117)  {\scriptsize $ A $}
     \put(56,155)  {\scriptsize $ A $}
     \put(38,117)  {\scriptsize $ A $}
     \put(10,41)  {\scriptsize $ A $}
     \put(33,192)  {\scriptsize $ U_i $}
     \put(33,38)  {\scriptsize $ U_j $}
     \put(6,137)  {\scriptsize $ M_{\kappa} $}
     \put(80,137)  {\scriptsize $ M_{\kappa} $}
     \put(89,19)  {\scriptsize $ e_j $}
     \put(89,203)  {\scriptsize $ r_i $}
    }\setlength{\unitlength}{1pt}}
  \end{picture}} 
\overset{(2)}{=} \sum_{i,j,\kappa} c_{i\kappa}
\raisebox{-91pt}{
\begin{picture}(85,190)
   \put(0,8){\scalebox{.75}{\includegraphics{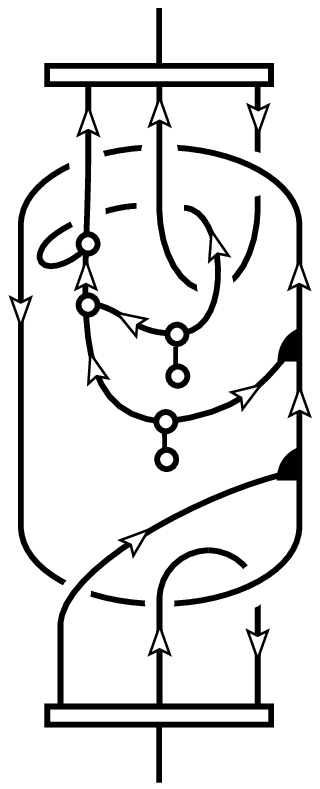}}}
   \put(0,8){
     \setlength{\unitlength}{.75pt}\put(-7,0){
     \put(48,-10)  {\scriptsize $ C_A $}
     \put(48,231)  {\scriptsize $ C_A $}
     \put(7,41)  {\scriptsize $ A $}
     \put(-12,138)  {\scriptsize $ M_{\kappa} $}
     \put(98,138)  {\scriptsize $ M_{\kappa} $}
     \put(98,107)  {\scriptsize $ M_{\kappa} $}
     \put(56,192)  {\scriptsize $ U_i $}
     \put(56,38)  {\scriptsize $ U_j $}
     \put(15,188)  {\scriptsize $ A $}
     \put(70,99)  {\scriptsize $ A $}
     \put(21,117)  {\scriptsize $ A $}
     \put(38,143)  {\scriptsize $ A $}
     \put(70,133)  {\scriptsize $ A $}
     \put(89,19)  {\scriptsize $ e_j $}
     \put(89,203)  {\scriptsize $ r_i $}
   }\setlength{\unitlength}{1pt}}
  \end{picture}}
\overset{(3)}{=} \sum_{i,j,\kappa} c_{i\kappa}
  \raisebox{-91pt}{
  \begin{picture}(70,190)
   \put(0,8){\scalebox{.75}{\includegraphics{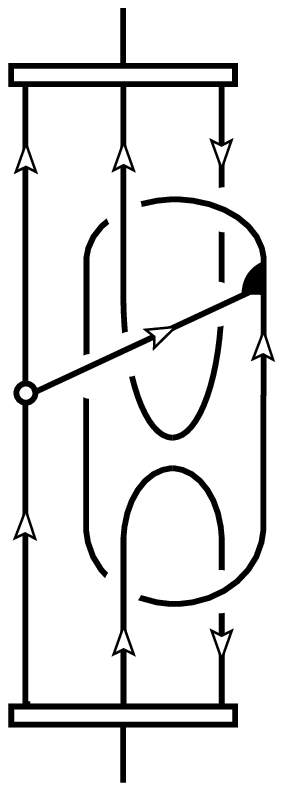}}}
   \put(0,8){
     \setlength{\unitlength}{.75pt}\put(-17,0){
     \put(48,-10)  {\scriptsize $ C_A $}
     \put(48,230)  {\scriptsize $ C_A $}
     \put(8,73)  {\scriptsize $ A $}
     \put(8,178)  {\scriptsize $ A $}
     \put(33,178)  {\scriptsize $ U_i $}
     \put(33,39)  {\scriptsize $ U_j $}
     \put(98,125)  {\scriptsize $ M_{\kappa} $}
     \put(89,19)  {\scriptsize $ e_j $}
     \put(89,203)  {\scriptsize $ r_i $}
    }\setlength{\unitlength}{1pt}}
  \end{picture}} 
$$
\be  
\overset{(4)}{=} \sum_{i,j,l,\kappa,\nu} 
  \frac{\dim(U_i) \dim(M_\kappa) \dim(U_l)}{\mathrm{Dim}(\Cc) \dim(U_j)}
\raisebox{-91pt}{
\begin{picture}(120,190)
   \put(0,8){\scalebox{.75}{\includegraphics{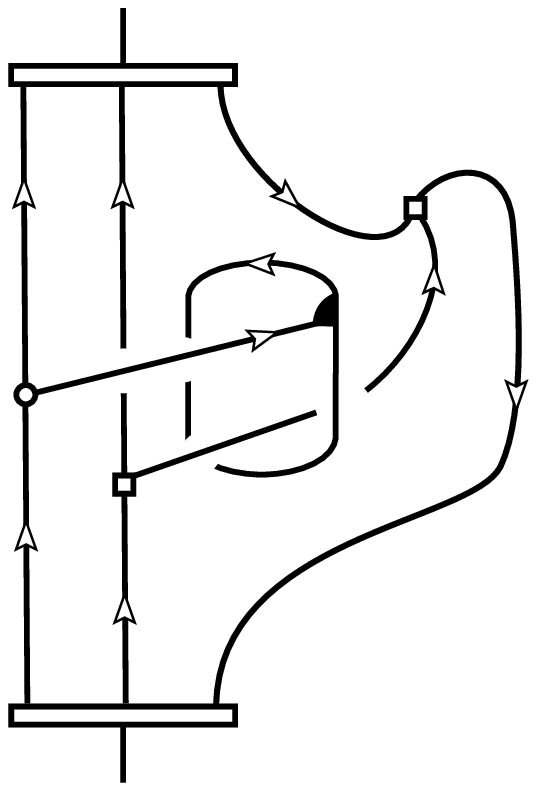}}}
   \put(0,8){
     \setlength{\unitlength}{.75pt}\put(-17,0){
     \put(48,-10)  {\scriptsize $ C_A $}
     \put(48,230)  {\scriptsize $ C_A $}
     \put(8,70)  {\scriptsize $ A $}
     \put(8,169)  {\scriptsize $ A $}
     \put(86,159)  {\scriptsize $ M_{\kappa} $}
     \put(33,169)  {\scriptsize $ U_i $}
     \put(96,182)  {\scriptsize $ U_i $}
     \put(33,49)  {\scriptsize $ U_j $}
     \put(146,111)  {\scriptsize $ U_j $}
     \put(146,142)  {\scriptsize $ U_l $}
     \put(90,118)  {\scriptsize $ A $}
     \put(38,85)  {\scriptsize $ \nu $}
     \put(122,168)  {\scriptsize $ \nu $}
     \put(89,19)  {\scriptsize $ e_j $}
     \put(89,203)  {\scriptsize $ r_i $}
   }\setlength{\unitlength}{1pt}}
  \end{picture}}
~.
\labl{eq:phi-inv-2-aux}
Equality (1) follows by substituting the definitions of $\varphi$ and $\bar\varphi$, and using $e_\kappa \cir r_\lambda = \delta_{\kappa,\lambda} P_{\otimes A}$. The $A$-ribbon in the idempotent $P_{\otimes A}$ can be rearranged to form the idempotent $P^l_A(U_i)$ using the representation property of $M_\kappa$ and that $A$ is symmetric Frobenius. This is done in step (2). In step (3) one uses that 
$r_i \cir (P^l_A(U_i) \oti \id_{U_i^\vee}) = r_i$, as well as the 
representation property of $A$ so that there is now only one $A$-ribbon
attached to the $M_\kappa$-ribbon. In step (4) the $U_i$ and 
$U_j$-ribbons are replaced by a sum over $U_l$ which amounts to the
decomposition of the tensor product $U_i^\vee \oti U_j$; the precise
identity employed is \cite[eqn.\,(4.8)]{defect} (or rather a
vertically reflected version thereof) for $\one$-$\one$-bimodules.
On the right hand side of \erf{eq:phi-inv-2-aux} one can now apply
lemma \ref{lem:M-loop}, and after cancelling all the 
factors and using that $\sum_{i \in \Ic} r_i \cir e_i = \id_{C_A}$
one arrives at the statement of the lemma.
\epf

\bigskip\noindent
{\it Proof of proposition \ref{prop:CA=TA}:}
\\[.3em]
Lemmas \ref{lem:phi-inv-1} and \ref{lem:phi-inv-2} imply that $\varphi$ is an isomorphism. It remains to check that it is an algebra map.
\\[.3em]
a) {\it $e_C$ is an algebra map:} 
Recall the definition of $e_C$ and $r_C$ in \erf{eq:ec-rc-def}.
By definition, $\eta_{C_A} = T(r_l \circ \eta_{R(A)})$. 
We have 
\be
e_C \circ \eta_{C_A} = T( \iota_l\circ r_l \circ \eta_{R(A)}) = 
T(\eta_{R(A)}) = \eta_{T(R(A))} ~,
\labl{eq:CA-TA-aux1}
where in the first step we used that $T$ is a functor, in the
second step we used \cite[lem.\,3.10]{corr} to
omit the idempotent $\iota_l\cir r_l$, and the third step is just the
definition of the unit of $T(R(A))$. 
For the multiplication we have, again
by definition, $m_{C_A}:= r_C \circ m_{T(R(A))} \circ (e_C\otimes e_C)$. 
Along the same lines as in \erf{eq:CA-TA-aux1} one computes
\be\bearll
e_{C} \circ m_{C_A} 
\overset{(1)}{=} T(\iota_l) \circ T(r_l) \circ 
T(m_{R(A)}) \circ \varphi_2 \circ (T(\iota_l) \otimes T(\iota_l))
\enl
\overset{(2)}{=} T(\iota_l) \circ T(r_l) \circ T(m_{R(A)}) \circ T(\iota_l \otimes \iota_l)
\circ \varphi_2 
\overset{(3)}{=} T\big(P_l(R(A))\circ m_{R(A)} \circ 
(\iota_l \otimes \iota_l)\big)
\circ \varphi_2 
\enl
\overset{(4)}{=} T(m_{R(A)} \circ (\iota_l\otimes \iota_l)) \circ \varphi_2
\overset{(5)}{=} m_{T(R(A))} \circ (e_C \otimes e_C).
\eear\ee
where in the first step the definitions in
\erf{eq:TB} and
\erf{eq:ec-rc-def} have been subsituted, and in
step (2) we used that $\varphi_2$ is a 
natural transformation, see section \ref{sec:functor-T}.
Step (4) is a consequence of \cite[lem.\,3.10]{corr}.
In step (5) one reverses step (2) and substitutes the definition of $e_C$.
\\[.3em]
b) $\varphi \circ m_{C_A} = m_{T_A} \circ (\varphi \otimes \varphi)$:
For $\varphi \circ m_{C_A}$ consider the equalities
$$
\varphi \circ m_{C_A} 
\overset{(1)}{=} \sum_{k,\kappa}  \,\,\,\,
  \raisebox{-66pt}{
  \begin{picture}(80,140)
   \put(0,8){\scalebox{.75}{\includegraphics{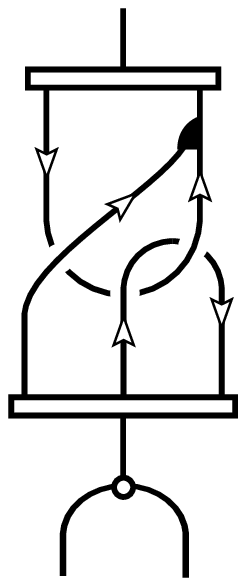}}}
   \put(0,8){
     \setlength{\unitlength}{.75pt}\put(-17,-105){
     \put(28,96)  {\scriptsize $ C_A $}
     \put(63,96)  {\scriptsize $ C_A $}
     \put(31,140)  {\scriptsize $ C_A $}
     \put(44,277)  {\scriptsize $ T_A $}
     \put(10,173)  {\scriptsize $ A $}
     \put(32,173)  {\scriptsize $ U_k $}
     \put(6,226)  {\scriptsize $ M_{\kappa} $}
     \put(89,156)  {\scriptsize $ e_k $}
     \put(84,249)  {\scriptsize $ r_{\kappa} $}
    }\setlength{\unitlength}{1pt}}
  \end{picture}} 
\overset{(2)}{=} \sum_{i,j,k,\kappa,\alpha}  
\raisebox{-86pt}{
\begin{picture}(110,180)
   \put(0,8){\scalebox{.75}{\includegraphics{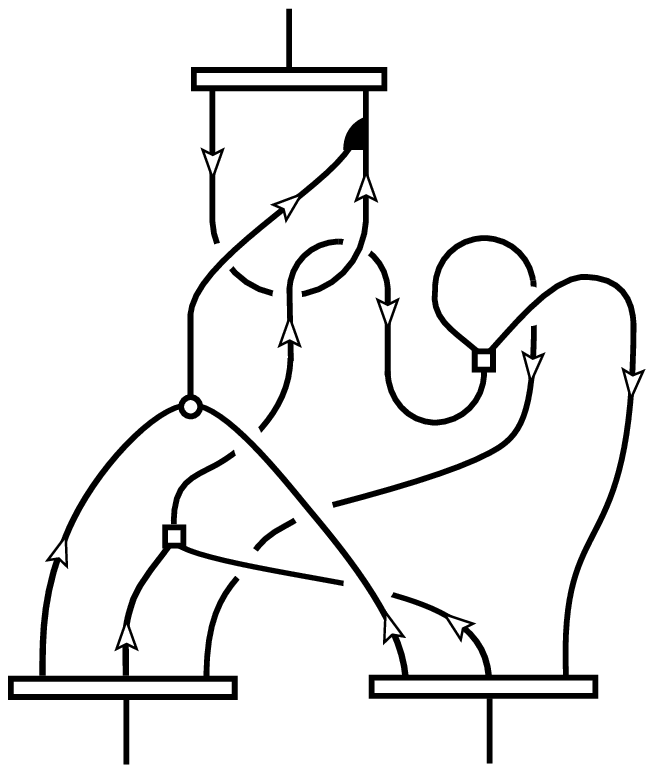}}}
   \put(0,8){
     \setlength{\unitlength}{.75pt}\put(-201,-129){
     \put(229,119)  {\scriptsize $ C_A $}
     \put(336,119)  {\scriptsize $ C_A $}
     \put(277,356)  {\scriptsize $ T_A $}
     \put(195,168)  {\scriptsize $ A $}
     \put(295,168)  {\scriptsize $ A $}
     \put(238,303)  {\scriptsize $ M_{\kappa} $}
     \put(216,165)  {\scriptsize $ U_i $}
     \put(360,245)  {\scriptsize $ U_i $}
     \put(343,165)  {\scriptsize $ U_j $}
     \put(390,240)  {\scriptsize $ U_j $}
     \put(262,251)  {\scriptsize $ U_k $}
     \put(238,255)  {\scriptsize $ A $}
     \put(249,185)  {\scriptsize $ \alpha $}
     \put(335,258)  {\scriptsize $ \alpha $}
     \put(271,150)  {\scriptsize $ e_i $}
     \put(376,150)  {\scriptsize $ e_j $}
     \put(315,326)  {\scriptsize $ r_{\kappa} $}
   }\setlength{\unitlength}{1pt}}
  \end{picture}}
$$
\be
\overset{(3)}{=} \sum_{i,j,\kappa}  \,\,\,\,
  \raisebox{-86pt}{
  \begin{picture}(150,180)
   \put(0,8){\scalebox{.75}{\includegraphics{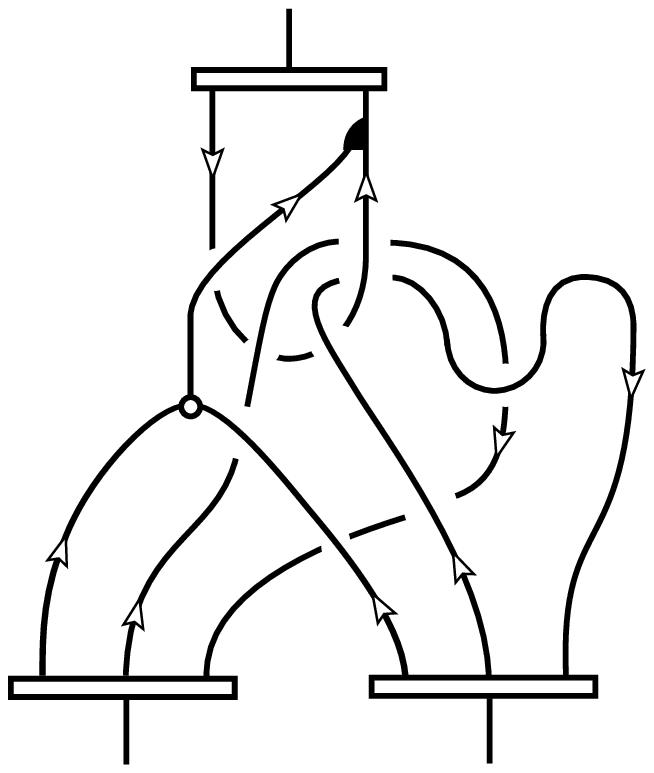}}}
   \put(0,8){
     \setlength{\unitlength}{.75pt}\put(-200,-129){
     \put(229,119)  {\scriptsize $ C_A $}
     \put(334,119)  {\scriptsize $ C_A $}
     \put(277,356)  {\scriptsize $ T_A $}
     \put(195,168)  {\scriptsize $ A $}
     \put(295,168)  {\scriptsize $ A $}
     \put(238,255)  {\scriptsize $ A $}
     \put(219,170)  {\scriptsize $ U_i $}
     \put(342,170)  {\scriptsize $ U_j $}
     \put(234,303)  {\scriptsize $ M_{\kappa} $}
     \put(271,150)  {\scriptsize $ e_i $}
     \put(376,150)  {\scriptsize $ e_j $}
     \put(315,326)  {\scriptsize $ r_{\kappa} $}
   }\setlength{\unitlength}{1pt}}
  \end{picture}}
~.
\labl{eq:CA=TA-aux1}
Here step (1) is the definition of $\varphi$. In step (2) we use part a) of the proof showing that $e_C$ is an algebra map, allowing us to 
replace the multiplication of $C_A$ by that of $T(R(A))$. In step (3)
the sum over $k$ and $\alpha$ is carried out, joining the two $U_i$-ribbons and the two $U_j$-ribbons, see e.g.\ \cite[eqn.\,(2.31)]{tft1}.
For $m_{T_A} \circ (\varphi \otimes \varphi)$ consider the equalities
$$
  m_{T_A} \circ (\varphi \otimes \varphi) \overset{(1)}{=} 
  \sum_{\mu, \sigma,\kappa, i,j}\,\,\,\,
  \raisebox{-86pt}{
  \begin{picture}(110,180)
   \put(0,8){\scalebox{.75}{\includegraphics{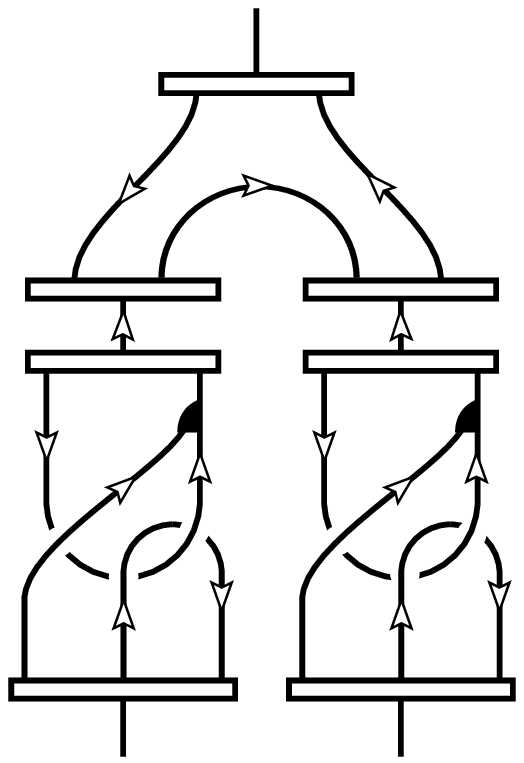}}}
   \put(0,8){
     \setlength{\unitlength}{.75pt}\put(-17,0){
     \put(45,-10)  {\scriptsize $ C_A $}
     \put(124,-10)  {\scriptsize $ C_A $}
     \put(80,224)  {\scriptsize $ T_A $}
     \put(10,39)  {\scriptsize $ A $}
     \put(90,39)  {\scriptsize $ A $}
     \put(32,39)  {\scriptsize $ U_i $}
     \put(112,39)  {\scriptsize $ U_j $}
     \put(7,90)  {\scriptsize $ M_{\kappa} $}
     \put(87,90)  {\scriptsize $ M_{\sigma} $}
     \put(60,123)  {\scriptsize $ T_A $}
     \put(140,123)  {\scriptsize $ T_A $}
     \put(28,166)  {\scriptsize $ M_{\mu} $}
     \put(135,166)  {\scriptsize $ M_{\mu} $}
     \put(84,175)  {\scriptsize $ M_{\mu} $}
     \put(7,19)  {\scriptsize $ e_i $}
     \put(168,19)  {\scriptsize $ e_j $} 
     \put(124,193)  {\scriptsize $ r_{\mu} $}
     }\setlength{\unitlength}{1pt}}
  \end{picture}} 
  \overset{(2)}{=}  \sum_{\kappa, i,j}
\raisebox{-86pt}{
\begin{picture}(110,180)
   \put(0,8){\scalebox{.75}{\includegraphics{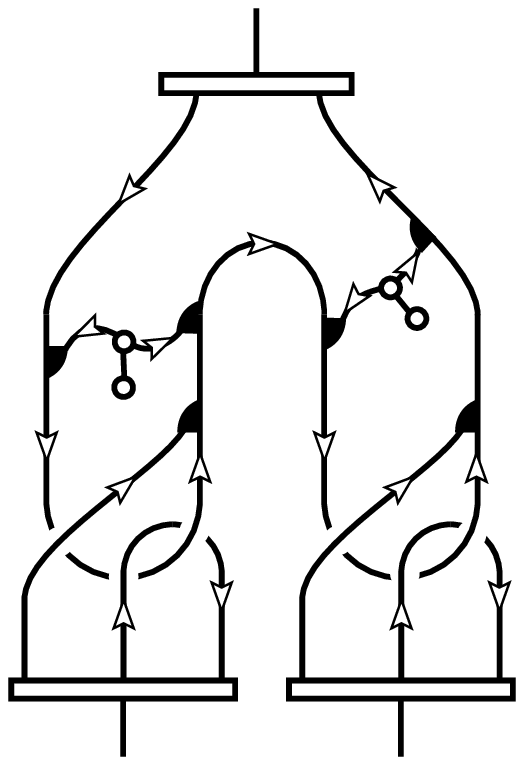}}}
   \put(0,8){
     \setlength{\unitlength}{.75pt}\put(-17,0){
     \put(45,-10)  {\scriptsize $ C_A $}
     \put(124,-10)  {\scriptsize $ C_A $}
     \put(80,224)  {\scriptsize $ T_A $}
     \put(10,39)  {\scriptsize $ A $}
     \put(90,39)  {\scriptsize $ A $}
     \put(32,39)  {\scriptsize $ U_i $}
     \put(112,39)  {\scriptsize $ U_j $}
     \put(28,166)  {\scriptsize $ M_{\kappa} $}
     \put(133,166)  {\scriptsize $ M_{\kappa} $}
     \put(54,128)  {\scriptsize $ A $}
     \put(112,142)  {\scriptsize $ A $}
     \put(123,146)  {\scriptsize $ A $}
     \put(6,19)  {\scriptsize $ e_i $}
     \put(168,19)  {\scriptsize $ e_j $} 
     \put(124,193)  {\scriptsize $ r_{\kappa} $}
     \put(37,133)  {\scriptsize $ \fbox{1} $}
     \put(75,96)  {\scriptsize $ \fbox{2} $}
     }\setlength{\unitlength}{1pt}}
  \end{picture}}
$$
\be
 \qquad \overset{(3)}{=} \sum_{\kappa,i,j}\,\,\,\,
  \raisebox{-76pt}{
  \begin{picture}(110,160)
   \put(0,8){\scalebox{.75}{\includegraphics{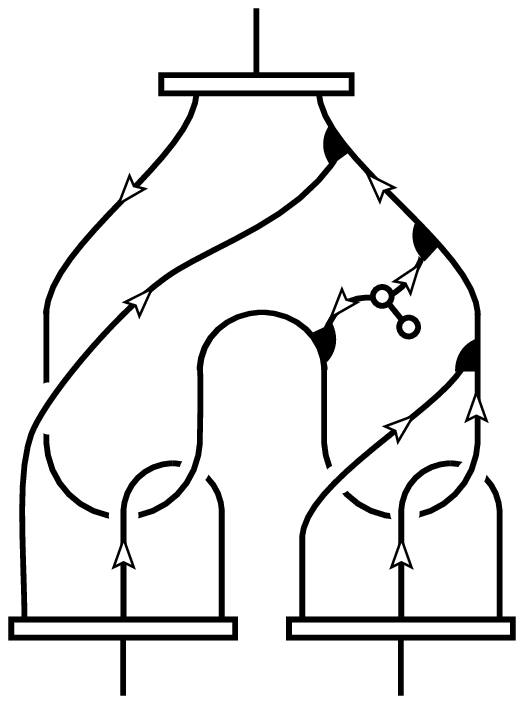}}}
   \put(0,8){
     \setlength{\unitlength}{.75pt}\put(-17,0){
     \put(45,-10)  {\scriptsize $ C_A $}
     \put(124,-10)  {\scriptsize $ C_A $}
     \put(80,207)  {\scriptsize $ T_A $}
     \put(10,39)  {\scriptsize $ A $}
     \put(90,39)  {\scriptsize $ A $}
     \put(123,127)  {\scriptsize $ A $}
     \put(101,119)  {\scriptsize $ A $}
     \put(32,39)  {\scriptsize $ U_i $}
     \put(112,39)  {\scriptsize $ U_j $}
     \put(29,150)  {\scriptsize $ M_{\kappa} $}
     \put(134,150)  {\scriptsize $ M_{\kappa} $}     
     \put(7,19)  {\scriptsize $ e_i $}
     \put(168,19)  {\scriptsize $ e_j $} 
     \put(124,176)  {\scriptsize $ r_{\kappa} $}
    }\setlength{\unitlength}{1pt}}
  \end{picture}} 
  \overset{(4)}{=}  \sum_{\kappa, i,j}
\raisebox{-76pt}{
\begin{picture}(110,160)
   \put(0,8){\scalebox{.75}{\includegraphics{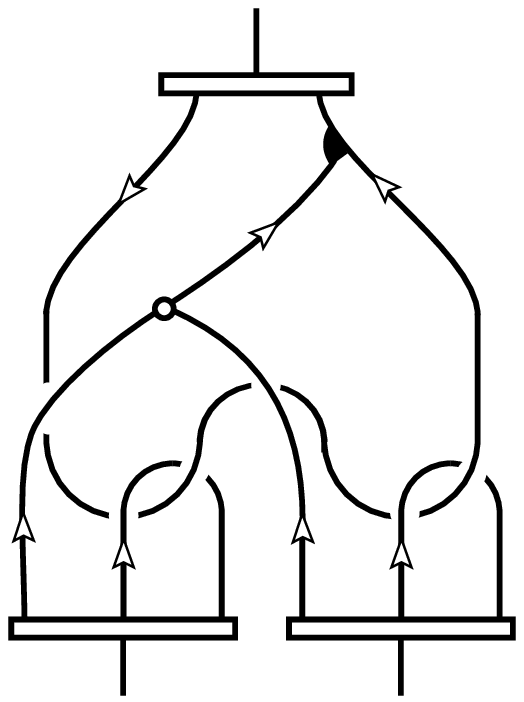}}}
   \put(0,8){
     \setlength{\unitlength}{.75pt}\put(-17,0){
     \put(45,-10)  {\scriptsize $ C_A $}
     \put(124,-10)  {\scriptsize $ C_A $}
     \put(80,207)  {\scriptsize $ T_A $}
     \put(6,45)  {\scriptsize $ A $}
     \put(86,45)  {\scriptsize $ A $}
     \put(29,150)  {\scriptsize $ M_{\kappa} $}
     \put(134,150)  {\scriptsize $ M_{\kappa} $}
     \put(32,39)  {\scriptsize $ U_i $}
     \put(112,39)  {\scriptsize $ U_j $}
     \put(76,136)  {\scriptsize $ A $}
     \put(7,19)  {\scriptsize $ e_i $}
     \put(168,19)  {\scriptsize $ e_j $} 
     \put(124,176)  {\scriptsize $ r_{\kappa} $}
   }\setlength{\unitlength}{1pt}}
  \end{picture}}
~.
\labl{eq:CA=TA-aux2}
The first equality follows again by substituting the definitions
and step (2) follows from 
$e_\mu \cir r_\kappa = \delta_{\mu,\kappa} P_{\otimes A}$.
In step (3) we have first removed the idempotent marked `1' by
rearranging it to become the idempotent $P^l_A(U_i)$ sitting on
top of the $e_i$ morphism, where it can be omitted. Then
the representation morphism marked `2' is dragged to the right,
and the representation property as well as that $A$ is symmetric
Frobenius is used to move the $A$-ribbon along the projector.
In step (4), one now
removes the remaining $P_{\otimes A}$ idempotent as before 
by rearranging it to be the idempotent $P^l_A(U_j)$ and
omitting it against the $e_j$ morphism. One also uses
once more the representation property of $M_\kappa$.
The result is easily
seen to agree with the right hand side of \erf{eq:CA=TA-aux1}.
\\[.3em]
c) $\varphi \circ \eta_{C_A} = \eta_{T_A}$:
This is an immediate consequence of combining
$e_i \cir \eta_{C_A} = \delta_{i,0}\,\eta_A \otimes \id_\one \otimes \id_\one$ with
the definition of $\varphi$ and using that
$\eta_{T_A} = \sum_\kappa r_\kappa \cir \tilde b_{M_\kappa}$.
\\[.3em]
Altogether we established that $\varphi$ (and hence also
$\bar\varphi$) is an
isomorphism of unital algebras.
\epf

\medskip

In the special case that $A=\one$, the fact that $\varphi$ is an algebra map already follows from the proof of theorem 5.19 in \cite{ko-cardy}.

\subsection{A surjection from $T(Z(A))$ to $A$}\label{sec:surj-TZA-A}

Let now $A$ be a haploid non-degenerate algebra in $\mathcal{C}$.
(Thus $A$ is in particular simple.)
Recall the definitions of $C_A$ in \erf{eq:CA-TA-def} 
and $e_i$, $r_i$ in \erf{eq:ei-ri-def_CA}. Define the morphisms
$\iota : C_A \rightarrow A$ and $\bar\iota : A \rightarrow C_A$
as
\be
  \iota = \sum_{i \in \Ic} (\id_A \oti \tilde d_{U_i}) \cir e_i
  \qquad \text{and} \qquad
  \bar\iota = \sum_{i \in \Ic} 
  \frac{\dim(A) \dim(U_i)}{\mathrm{Dim}(\Cc)} \,
  r_i \cir (\id_A \oti b_{U_i}) ~.
\labl{eq:iota-A_CA-def}
As mentioned in section \ref{sec:pre-nondeg}, for simple 
non-degenerate $A$ we automatically have $\dim(A) \neq 0$.

\begin{lemma}  \label{lemma:iota-iota-id}
$\iota \cir \bar{\iota} = \id_A$.
\end{lemma}
\pf 
Let $\{\,x_\alpha^i\,\}$ be a basis of 
$\Hom(U_i,A)$
and let $\{\,\bar x_\alpha^i\,\}$ be the basis 
of $\Hom(A,U_i)$ dual to $x_\alpha^i$ in the sense that
$\bar x_\alpha^i \cir x_\beta^i = \delta_{\alpha,\beta} \, \id_{U_i}$.
Since $A$ is haploid, for
$i=0$ there is only one basis vector in each space,
and we choose $x_1^0 = \eta_A$ and
$\bar x_1^0 = \dim(A)^{-1} \eps_{A}$.
We have
\bea
\iota \circ \bar{\iota} = \sum_{i} 
   \frac{\dim(A) \dim(U_i)}{\mathrm{Dim}(\Cc)}\,\,\,
  \raisebox{-34pt}{
  \begin{picture}(60,75)
   \put(0,8){\scalebox{.75}{\includegraphics{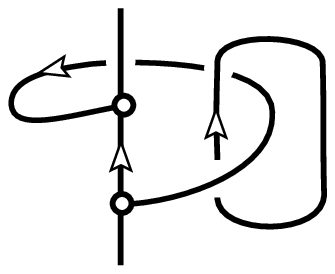}}}
   \put(0,8){
     \setlength{\unitlength}{.75pt}\put(-218,-429){
     \put(245,421)  {\scriptsize $ A $}
     \put(235,458)  {\scriptsize $ A $}
     \put(261,468)  {\scriptsize $ U_i $}
     \put(245,510)  {\scriptsize $ A $}
     }\setlength{\unitlength}{1pt}}
  \end{picture}} 
\enl
=  \frac{\dim(A)}{\mathrm{Dim}(\Cc)}\sum_{i,k,\alpha} \dim(U_i) 
\raisebox{-51pt}{
\begin{picture}(70,110)
   \put(0,8){\scalebox{.75}{\includegraphics{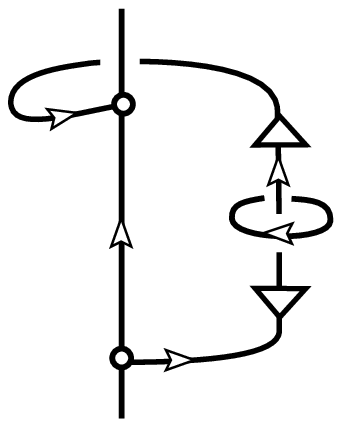}}}
   \put(0,8){
     \setlength{\unitlength}{.75pt}\put(-102,-187){
     \put(129,178)  {\scriptsize $ A $}
     \put(129,313)  {\scriptsize $ A $}
     \put(119,240)  {\scriptsize $ A $}
     \put(146,215)  {\scriptsize $ A $}
     \put(161,297)  {\scriptsize $ A $}
     \put(150,245)  {\scriptsize $ U_i $}
     \put(187,256)  {\scriptsize $ U_k $}
     \put(159,270)  {\scriptsize $ x^k_\alpha $}
     \put(159,221)  {\scriptsize $ \bar x^k_{\alpha} $}
     }\setlength{\unitlength}{1pt}}
  \end{picture}}
~= \dim(A)\sum_{k,\alpha} \delta_{k,0}~~~~
\raisebox{-34pt}{
  \begin{picture}(60,75)
   \put(0,8){\scalebox{.75}{\includegraphics{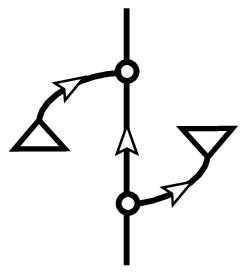}}}
   \put(0,8){
     \setlength{\unitlength}{.75pt}\put(-218,-429){
     \put(245,421)  {\scriptsize $ A $}
     \put(226,490)  {\scriptsize $ A $}
     \put(267,442)  {\scriptsize $ A $}
     \put(245,510)  {\scriptsize $ A $}
     \put(203,465)  {\scriptsize $ x^0_\alpha $}
     \put(286,464)  {\scriptsize $ \bar x^0_{\alpha} $}
     }\setlength{\unitlength}{1pt}}
  \end{picture}} 
  ~~.
\eear\ee
Since $A$ is haploid, the last sum over $\alpha$ only
contains one term, and by our convention on
$x_1^0$ and $\bar x_1^0$ the right hand side is then
just equal to $\id_A$.
\epf

\begin{lemma}  \label{lemma:iota-algebra-hom}
$\iota$ is an algebra map.
\end{lemma}
\pf 
We have
\bea
    \iota \cir m_{C_A} 
     \overset{(1)}{=} \sum_{i,j,k,\alpha} \,\,\,\,
 \raisebox{-52pt}{
  \begin{picture}(120,112)
   \put(0,8){\scalebox{.75}{\includegraphics{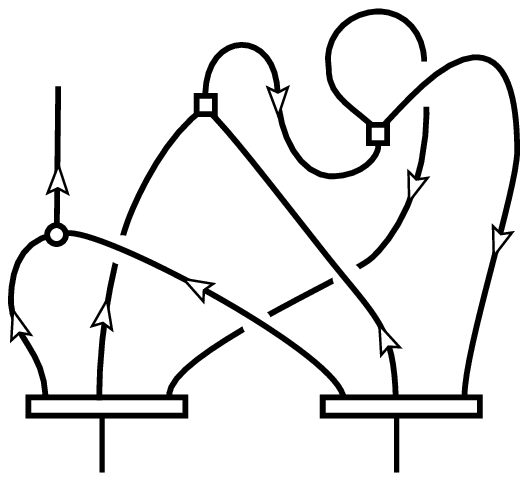}}}
   \put(0,8){
     \setlength{\unitlength}{.75pt}\put(-97,-62){
     \put(115,52)  {\scriptsize $ C_A $}
     \put(202,52)  {\scriptsize $ C_A $}
     \put(85,103,)  {\scriptsize $ A $}
     \put(155,125)  {\scriptsize $ A $}
     \put(107,181)  {\scriptsize $ A $}
     \put(110,114)  {\scriptsize $ U_i $}
     \put(223,144)  {\scriptsize $ U_i $}
     \put(215,99)  {\scriptsize $ U_j $}
     \put(250,142)  {\scriptsize $ U_j $}
     \put(161,192)  {\scriptsize $ U_k $}
     \put(150,155)  {\scriptsize $ \alpha $}
     \put(200,173)  {\scriptsize $ \alpha $}
     \put(89,81)  {\scriptsize $ e_i $}
     \put(240,81)  {\scriptsize $ e_j $}
  }\setlength{\unitlength}{1pt}}
  \end{picture}}
 \overset{(2)}{=} \sum_{i,j} \,\,\,\,
  \raisebox{-41pt}{
  \begin{picture}(95,90)
   \put(0,8){\scalebox{.75}{\includegraphics{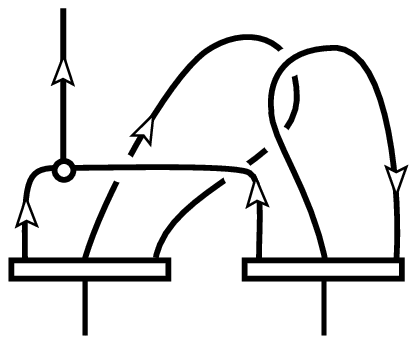}}}
   \put(0,8){
     \setlength{\unitlength}{.75pt}\put(-97,-62){
     \put(115,52)  {\scriptsize $ C_A $}
     \put(202,52)  {\scriptsize $ C_A $}
     \put(88,99)  {\scriptsize $ A $}
     \put(155,95)  {\scriptsize $ A $}
     \put(107,161)  {\scriptsize $ A $}
     \put(120,124)  {\scriptsize $ U_i $}
     \put(217,107)  {\scriptsize $ U_j $}
     \put(87,81)  {\scriptsize $ e_i $}
     \put(217,81)  {\scriptsize $ e_j $}
    }\setlength{\unitlength}{1pt}}
  \end{picture}}
\enl
\overset{(3)}{=} 
m_A \cir (\iota \oti \iota)
~.
\eear\ee
In the first step, as in the second step of \erf{eq:CA=TA-aux1}
we use that $e_C$ is an algebra map to replace the multiplication
of $C_A$ by that of $T(R(A))$. In step (2) we carry out the sum
over $k,\alpha$ (as in step (3) of \erf{eq:CA=TA-aux1}),
and equality (3) is then immediate by deforming the ribbons.
For the unit we get, using also \erf{eq:CA-TA-aux1},
\be
  \iota \cir \eta_{C_A} = 
  \sum_{i \in \Ic} (\id_A \oti \tilde d_{U_i}) \cir \eta_{T(R(A))}
  = \eta_A
  ~.
\ee
\epf

\subsection{Haploid representatives of Morita classes}\label{sec:hapl-rep}

The following proposition establishes that the Morita class of a simple non-degenerate algebra always contains a haploid representative. This fact will be used in the proof of theorem \ref{thm:main}.

\begin{prop}\label{prop:MxM-prop2}
Let $A$ be a simple non-degenerate algebra in a modular
tensor category.
\\[.3em]
(i)\phantom{i} 
Given a left $A$-module $M$ with $\dim(M) \neq 0$, 
the algebra $M^{\vee} \otimes_A M$ is simple, non-degenerate,
and Morita-equivalent to $A$.
\\[.3em]
(ii) $A$ is Morita-equivalent to a haploid non-degenerate algebra.
\end{prop}
\pf 
(i) The algebra structure on $B := M^{\vee} \otimes_A M$ was given
in lemma \ref{lem:MxM-prop1}\,(i). As mentioned in section
\ref{sec:pre-nondeg} for simple non-degenerate $A$ we 
automatically have $\dim(A) \neq 0$. Thus by 
lemma \ref{lem:nondeg_vs_Frob}\,(i,iv) $A$ is simple 
normalised-special
symmetric Frobenius, and we can apply \cite[prop.\,2.13]{tft2} to
conclude that also $B$ is simple normalised-special
symmetric Frobenius (this uses $\dim(M) \neq 0$). 
By lemma \ref{lem:nondeg_vs_Frob}\,(ii)
$B$ is then in particular simple and non-degenerate.
That $A$ and $B$ are Morita equivalent follows from
\cite[thm.\,2.14]{tft2}.
\\[.3em]
(ii) Let $M$ be a simple left $A$-module. Applying 
\cite[lem.\,2.6]{defect} in the special case of $A$-$\one$-bimodules
shows that $\dim(M) \neq 0$. By part (i), $M^{\vee} \otimes_A M$
is Morita-equivalent to $A$ and by lemma \ref{lem:MxM-prop1}\,(ii),
$M^{\vee} \otimes_A M$ is haploid.
\epf

\medskip

The proposition essentially also follows from \cite[sect.\,3.3]{ost},
which however works in a slightly different setting. Note also that
the above proof does not make use of the modularity (or even the braiding)
of $\Cc$. We restrict our attention to the modular case because we
want to avoid changing the categorial framework repeatedly.

\medskip

We have now gathered all the ingredients to complete the proof 
of theorem \ref{thm:main}.

\bigskip\noindent
{\it Proof of (ii)$\,\Rightarrow$(i) in theorem \ref{thm:main}}:
\\[.3em]
We are given two simple non-degenerate algebras $A$, $B$ in
$\Cc$ such that $Z(A) \cong Z(B)$ as algebras. By proposition
\ref{prop:MxM-prop2}\,(ii) we can find a haploid non-degenerate 
algebra $B'$ that is Morita equivalent to $B$. 
To prove that $A$ and $B$ are Morita-equivalent it is enough to
show that $A$ and $B'$ are Morita-equivalent. 
In section \ref{sec:i-to-ii} we have 
established that (i)$\Rightarrow$(ii)
in theorem \ref{thm:main}, and so $Z(B) \cong Z(B')$ as algebras.
Without loss of generality we can thus assume that $B$ is haploid.
\\[.3em]
{\it a) A surjective algebra map from $T_A$ to $B$:}
Let $f: Z(A) \rightarrow Z(B)$ be an algebra isomorphism.
We define a map $h: T_A \rightarrow B$ by the following 
composition of maps
\be
h = T_A \xrightarrow{\,\,\, \bar{\varphi} \,\,\,} 
T(Z(A)) \xrightarrow{~T(f)~} T(Z(B)) \xrightarrow{\,\,\, 
\iota \,\,\, } B ~~,
\ee
where $\bar\varphi$ was defined in \erf{eq:phi-phi-def} and
$\iota$ in \erf{eq:iota-A_CA-def}. 
By (the proof of) proposition \ref{prop:CA=TA}, $\bar\varphi$
is an algebra map, $T(f)$ is an algebra map since $T$ is 
a tensor functor, and $\iota$ is an algebra map according
to lemma \ref{lemma:iota-algebra-hom}. Thus $h$ is an algebra map.
Let $\bar{h}:= \varphi \circ T(f^{-1}) \circ \bar{\iota}$. 
Then by lemmas \ref{lemma:iota-iota-id} and \ref{lem:phi-inv-2}, 
we obtain $h \circ \bar{h} = \id_B$. Thus $h$ is also surjective.

Let 
$j_\kappa : M_{\kappa}^{\vee} \otimes_A M_{\kappa} \rightarrow T_A$ and 
$\pi_\kappa : T_A \rightarrow M_{\kappa}^{\vee} \otimes_A M_{\kappa}$
the embedding and projection for the subobject 
$M_{\kappa}^{\vee} \otimes_A M_{\kappa}$ of $T_A$.
Define $S \subset \Jc$ to consist of all $\kappa$ such that
$h \cir j_\kappa \neq 0$, and set
$T_A' = \bigoplus_{\kappa \in S} M_{\kappa}^{\vee} \otimes_A M_{\kappa}$.
Let $j' = \bigoplus_{\kappa \in S} j_\kappa : T_A' \rightarrow T_A$ 
be the embedding of the subobject $T_A'$ 
into $T_A$ and $\pi' : T_A \rightarrow T_A'$ the projection onto $T_A'$.
Let $h' = h \cir j'$, i.e.\ $h'$ is the restriction of $h$ to $T_A'$. 
\\[.3em]
{\it b) $h'$ is an algebra map:} 
Note that $j'$ obeys $j' \cir m_{T_A'} = m_{T_A} \cir (j' \oti j')$.
(However, for $S \neq \Jc$ $j'$ does not preserve the unit.)
Since $h$ is an algebra map it follows that also
$h' \cir m_{T_A'} = m_{B} \cir (h' \oti h')$. It remains to verify
that $h'$ preserves the unit. 
Note that $\eta_{T'_A} = \pi' \cir \eta_{T_A}$ and hence
\be\bearll
  h' \cir \eta_{T_A'} 
  \etb= h \cir j' \cir \pi' \cir \eta_{T_A}
  = \sum_{\kappa \in S} h \cir j_\kappa \cir \pi_\kappa \cir \eta_{T_A}
  = \sum_{\kappa \in \Jc} h \cir j_\kappa \cir \pi_\kappa \cir \eta_{T_A}
\\
  \etb= h \cir \eta_{T_A} = \eta_B \ .
\eear\ee
\\[.3em]
{\it c) $h'$ is surjective:} Suppose that $f \cir h' = 0$ for some morphism
$f : B \rightarrow U$ and some object $U$. Then
\be
  f \cir h = 
  \sum_{\kappa \in \Jc} f \cir h \cir j_\kappa \cir \pi_\kappa
  = \sum_{\kappa \in S} f \cir h \cir j_\kappa \cir \pi_\kappa
  = f \cir h \cir j' \cir \pi'
  = f \cir h' \cir \pi' = 0 ~.
\ee
Since $h$ is surjective, this implies that $f=0$. Altogether we see
that $f\cir h' = 0 \,\Rightarrow\, f\,{=}\,0$ and thus also $h'$ is
surjective.
\\[.3em]
{\it d) $h'$ is injective:}
Denote by $m_\kappa$ and $\eta_\kappa$ the multiplication and unit
of $M^\vee_\kappa \otimes_A M_\kappa$.
Just as was the case for $j'$, the morphism $j_\kappa$
obeys $j_\kappa \cir m_{\kappa} = m_{T_A} \cir (j_\kappa \oti j_\kappa)$.
This implies that the kernel of $j_\kappa$ will be a sub-bimodule
of $M^\vee_\kappa \otimes_A M_\kappa$, seen as a bimodule over itself.
The same holds for the combination $h' \cir j_\kappa$.
But $M^\vee_\kappa \otimes_A M_\kappa$ is simple, and hence $h' \cir j_\kappa$
is either injective or zero. In particular, for $\kappa\in S$, $h' \cir j_\kappa$
is injective. 

By assumption, $B$ is haploid and there exist constants 
$\lambda_\kappa \in \Cb$
such that $h' \cir j_\kappa \cir \eta_\kappa = \lambda_\kappa \eta_B$. 
Let $U$ be an object in $\Cc$ and $f : U \rightarrow T_A'$ a morphism. Suppose that $h' \cir f = 0$. Then
\bea
h' \cir f = 0
~\overset{(1)}\Rightarrow~
m_B \cir (\lambda_\kappa \eta_B \oti \id_B) \cir h' \cir f = 0
~\overset{(2)}\Rightarrow~
m_B \cir ((h' \cir j_\kappa \cir \eta_\kappa) \oti h') \cir f = 0
\enl
\overset{(3)}\Rightarrow~
h' \cir m_{T_A'} \cir ((j_\kappa \cir \eta_\kappa) \oti \id_{T_A'}) \cir f = 0
~\overset{(4)}\Rightarrow~
h' \cir j_\kappa \cir m_\kappa \cir (\eta_\kappa \oti \pi_\kappa) \cir f = 0 ~.
\enl
\overset{(5)}\Rightarrow~
h' \cir j_\kappa \cir \pi_\kappa \cir f = 0
~\overset{(6)}\Rightarrow~
\pi_\kappa \cir f = 0 ~\text{for all}~ \kappa \In S
~\overset{(7)}\Rightarrow~
\sum_{\kappa\in S} j_\kappa \cir \pi_\kappa \cir f = 0
\enl
\overset{(8)}\Rightarrow~
\id_{T_A'} \cir f = 0
\eear\ee
Step (1) follows from the unit property of $B$, in step (2) the above observation on the relation between $\eta_B$ and $\eta_\kappa$ is substituted, and step (3) follows since $h'$ is an algebra map. To see implication (4) one observes that 
$m_{T_A'} \cir (j_\kappa  \oti \id_{T_A'})
= j_\kappa \cir m_\kappa \cir (\id \oti \pi_\kappa)$,
step (5) is the unit property of $M^\vee_\kappa \otimes_A M_\kappa$, and
step (6) is implied by injectivity of $h' \cir j_\kappa$. Steps (7) and (8)
are clear.
Altogether, $h' \cir f = 0$ implies $f=0$, and hence $h'$ is injective. 
\\[.3em]
{\it e) $A$ and $B$ are Morita-equivalent:}
Combining parts b), c) and d) we see that $h' : T_A' \rightarrow B$ 
is a bijection of algebras. Since $B$ is haploid, $T_A'$ can only consist
of one summand, i.e.\ $|S|=1$. Let $\kappa$ be the unique element of $S$.
Then $h'$ is a bijection of algebras between
$M^\vee_\kappa \otimes_A M_\kappa$ and $B$. 
By proposition \ref{prop:MxM-prop2}\,(i), the algebra 
$M_{\kappa}^{\vee} \otimes_A M_{\kappa}$ is Morita equivalent
to $A$ and thus also $B$ is Morita equivalent to $A$.
\epf

\bigskip
\noindent
{\bf Acknowledgements:}\\
We would like to thank the `Zentrum f\"ur Mathematische Physik'
in Hamburg for organising the workshop
on the Geometric Langlands Program in July, where this work 
was started.
We would also like to thank 
Yi-Zhi Huang,
Michael M\"uger,
J\"urgen Fuchs, and especially
Christoph Schweigert
for helpful discussions and comments on the draft.
The research of IR was partially supported by
the EPSRC First Grant EP/E005047/1, the PPARC rolling grant
PP/C507145/1 and the Marie Curie network `Superstring Theory'
(MRTN-CT-2004-512194).




\begin{thebibliography}{XXX}

\bibitem[AF]{and-ful} F.W.~Anderson and K.R.~Fuller,
{\it Rings and Categories of Modules},
GTM 13, Springer 1973.

\bibitem[B]{bichon} J.~Bichon,
{\it Cosovereign Hopf algebras},
J.\ Pure Appl.\ Alg.\ {\bf 157} (2001) 121--133 [math.QA/9902030].

\bibitem[BK]{baki} B.~Bakalov and A.A.~Kirillov, 
{\it Lectures on Tensor Categories and Modular Functors} 
(American Mathematical Society, Providence 2001).

\bibitem[ENO]{eno} P.I.~Etingof, D.~Nikshych, and V.~Ostrik, 
{\it On fusion categories}, 
Ann.\ Math.\ {\bf 162} (2005) 581--642 [math.QA/0203060].
 
\bibitem[Fj1]{tft5} J.~Fjelstad, J.~Fuchs, I.~Runkel and C.~Schweigert,
{\it TFT construction of RCFT correlators. V: Proof of modular
invariance  and  factorisation},
Theo.\ and Appl.\ of\ Cat.\ {\bf 16} (2006) 342--433
[hep-th/0503194].

\bibitem[Fj2]{tftcft} J.~Fjelstad, J.~Fuchs, I.~Runkel and C.~Schweigert,
{\it Topological and conformal field theory as Frobenius algebras},
Contemp.\ Math.\ {\bf 431} (2007) 225--247 [math.CT/0512076].

\bibitem[Fj3]{unique} J.~Fjelstad, J.~Fuchs, I.~Runkel and C.~Schweigert,
{\it Uniqueness of open/closed rational CFT with given algebra of open states},
hep-th/0612306.
  
\bibitem[Fr1]{corr} J.~Fr\"ohlich, J.~Fuchs, I.~Runkel and C.~Schweigert,
{\it Correspondences of ribbon categories},
Adv.\ Math.\ {\bf 199} (2006) 192--329 [math.CT/0309465].

\bibitem[Fr2]{defect} J.~Fr\"ohlich, J.~Fuchs, I.~Runkel and C.~Schweigert,
{\it Duality and defects in rational conformal field theory},
Nucl.\ Phys.\  B {\bf 763} (2007) 354--430 [hep-th/0607247].
  
\bibitem[Fu1]{tft1} J.~Fuchs, I.~Runkel and C.~Schweigert,
{\it TFT construction of RCFT correlators. I: Partition functions},
Nucl.\ Phys.\ B {\bf 646} (2002) 353--497 [hep-th/0204148].

\bibitem[Fu2]{tft2} J.~Fuchs, I.~Runkel and C.~Schweigert,
{\it TFT construction of RCFT correlators. II: Unoriented world sheets},
Nucl.\ Phys.\  B {\bf 678} (2004) 511--637
[hep-th/0306164].
  
\bibitem[Fu3]{bim} J.~Fuchs, I.~Runkel and C.~Schweigert,
{\it The fusion algebra of bimodule categories},
Appl.\ Cat.\ Str.\ {\bf 16} (2008) 123--140 
[math.CT/0701223].

\bibitem[FS]{fs-cat} J.~Fuchs and C.~Schweigert,
{\it Category theory for conformal boundary conditions},
Fields Institute Commun.\ {\bf 39} (2003) 25--70 [math.CT/0106050].

\bibitem[H1]{hu1} Y.-Z.~Huang,
{\it Vertex operator algebras and the Verlinde conjecture},
Commun. Contemp. Math. {\bf 10} (2008) 103--154 [math.QA/0406291].
  
\bibitem[H2]{hu2} Y.-Z.~Huang,
{\it Rigidity and modularity of vertex tensor categories},
math.QA/0502533.

\bibitem[HK1]{hk1} Y.-Z.~Huang and L.~Kong,
{\it Open-string vertex algebras, tensor categories and operads},
Commun.\ Math.\ Phys.\ {\bf 250} (2004) 433--471 [math.QA/0308248].

\bibitem[HK2]{hk2} Y.-Z.~Huang and L.~Kong,
{\it Full field algebras},
Commun.\ Math.\ Phys.\  {\bf 272} (2007) 345--396 [math.QA/0511328].
   
\bibitem[JS]{js} A.~Joyal and R.~Street, 
{\it Braided tensor categories}, 
Adv.\ Math.\ {\bf 102} (1993) 20--78.

\bibitem[Kg1]{ko-ffa} L.~Kong,
{\it Full field algebras, operads and tensor  categories}, 
Adv.\ Math.\ {\bf 213} (2007) no. 1, 271--340 [math.QA/0603065]

\bibitem[Kg2]{ko-ocfa} L.~Kong, 
{\it Open-closed field algebras}, 
Commun.\ Math. \ Phys. {\bf 280} 207--261 [math.QA/0610293]. 

\bibitem[Kg3]{ko-cardy} L.~Kong,
{\it Cardy condition for open-closed field algebras},
math.QA/0612255.

\bibitem[Kk]{kock} J.~Kock, 
{\it Frobenius Algebras and 2D Topological Quantum 
  Field Theories}, Cambridge University Press, Cambridge 2003.

\bibitem[KLM]{klm} Y.~Kawahigashi, R.~Longo and M.~M\"uger,
{\it Multi-interval subfactors and modularity of representations 
in conformal field theory},
Commun.\ Math.\ Phys.\ {\bf 219} (2001) 631--669 [math.OA/9903104].
  
\bibitem[KO]{kios} A.J.~Kirillov and V.~Ostrik,
{\it On q-analog of McKay correspondence and ADE classification of 
$\widehat sl(2)$ conformal field theories},
Adv.\ Math.\ {\bf 171} (2002) 183--227 [math.QA/0101219].

\bibitem[KR]{us2} L.~Kong and I.~Runkel, {\it Cardy algebras and 
sewing constrains, I}, 0807.3356 [math.QA]. 

\bibitem[LP]{lp} A.D.~Lauda and H.~Pfeiffer,
{\it Open-closed strings: Two-dimensional extended TQFTs 
and Frobenius algebras}, Topology Appl.\ {\bf 155} (2008) 623--666. 
[math/0510664].

\bibitem[LR]{lr} R.~Longo and K.H.~Rehren,
{\it Local fields in boundary conformal QFT},
Rev.\ Math.\ Phys.\  {\bf 16} (2004) 909--960 [math-ph/0405067].
  
\bibitem[Mo]{mo} G.W.~Moore,
{\it Some comments on branes, G-flux, and K-theory},
Int.\ J.\ Mod.\ Phys.\  A {\bf 16} (2001) 936--944
[hep-th/0012007].
  
\bibitem[MS]{mose} G.W.~Moore and G.~Segal,
{\it D-branes and K-theory in 2D topological field theory},
hep-th/0609042.
  
\bibitem[M\"u1]{mug2} M.~M\"uger,
{\it From Subfactors to Categories and Topology II. 
The quantum double of tensor categories and subfactors},
J.\ Pure Appl.\ Alg.\ {\bf 180} (2003) 159--219 [math.CT/0111205].

\bibitem[M\"u2]{mug-conf} M.~M\"uger,
talk at workshop `Quantum Structures' (Leipzig, 28.\ June 2007),
preprint in preparation.

\bibitem[O]{ost} V.~Ostrik, 
{\it Module categories, weak Hopf algebras and modular invariants}, 
Transform.\ Groups {\bf 8} (2003) 177--206 [math.QA/0111139].

\bibitem[SS]{ss} A.N.~Schellekens and Y.S.~Stanev,
{\it Trace formulas for annuli},
JHEP {\bf 0112} (2001) 012 [hep-th/0108035].
  
\bibitem[T]{tur} V.G.~Turaev, 
{\it Quantum Invariants of Knots and 3-Manifolds}
(de Gruyter, New York 1994).

\bibitem[VZ]{vz} F.~Van Oystaeyen and Y.H.~Zhang, 
{\it The Brauer group of a braided monoidal category}, 
J.\ Algebra {\bf 202} (1998) 96--128.

\end{thebibliography}
\end{document}